\documentclass[12pt]{article}
\usepackage{amsfonts}
\usepackage{graphicx}
\def\Z{{\mathbb Z}}
\def\Q{{\mathbb Q}}
\def\mI{{\mathcal I}}
\def\mJ{{\mathcal J}}
\def\mL{{\mathcal L}}
\def\mU{{\mathcal U}}
\def\mG{{\mathcal G}}
\def\mod{\hbox{ mod }}
\def\product{\prod}
\def\pf{\par{\bf Proof:  }}
\def\pfof#1{\par{\bf Proof of #1:  }}
\def\eop{\hfill $\square$}
\newtheorem{theorem}{Theorem}
\newtheorem{lemma}{Lemma}
\newtheorem{slemma}{Lemma}
\begin{document}

\title{\Huge The Cube Recurrence}
\author{\LARGE Gabriel D. Carroll \normalsize \\ Harvard University \\ Cambridge, Massachusetts 02138 \\ \texttt{gcarroll@fas.harvard.edu} \\ \\ \LARGE David Speyer \normalsize \\ Department of Mathematics \\ University of California at Berkeley \\ Berkeley, California 94720 \\ \texttt{speyer@math.berkeley.edu}  \\}
\date{\small Submitted: December 31, 2002 \small MR Subject Classifications: 05A15, 05E99, 11B83}
\maketitle

\begin{abstract}
We construct a combinatorial model that is described by the cube recurrence, a nonlinear recurrence relation introduced by Propp, which generates families of Laurent polynomials indexed by points in $\Z^3$.  In the process, we prove several conjectures of Propp and of Fomin and Zelevinsky, and we obtain a combinatorial interpretation for the terms of Gale-Robinson sequences.  We also indicate how the model might be used to obtain some interesting results about perfect matchings of certain bipartite planar graphs.
\end{abstract}

\section{Introduction}
Consider a family of rational functions $f_{i,j,k}$, indexed by $(i,j,k) \in \Z^3$ with $k \geq -1$, and given by the initial conditions $f_{i,j,k} = x_{i,j,k}$ (a formal variable) for $k = -1, 0$, and $$f_{i,j,k-1}f_{i,j,k+1} = f_{i-1,j,k}f_{i+1,j,k} + f_{i,j-1,k}f_{i,j+1,k}.\quad (k \geq 0)$$  This is the {\em octahedron recurrence}, which has connections with the Hirota equation in physics, with Dodgson's condensation method of evaluating determinants, with alter\-nating-sign matrices, and with domino tilings of Aztec diamonds (see \cite{Zabrodin}, \cite{Proppfaces}, \cite{RobRum}, \cite{EKLP}, respectively).  It turns out that every $f_{i,j,k}$ is a Laurent polynomial in the initial $x_{i,j,k}$, i.e. a polynomial in the variables $x_{i,j,k}, x_{i,j,k}^{-1}$.  Sergey Fomin and Andrei Zelevinsky, using techniques from the theory of cluster algebras, proved in \cite{FomZel} that the recurrence again generates Laurent polynomials for a large variety of other sets of initial conditions (i.e. sets of points $(i,j,k)$ for which we designate $f_{i,j,k} = x_{i,j,k}$). In \cite{Speyer}, David Speyer showed further that all such polynomials could be interpreted as enumerating perfect matchings of suitable bipartite planar graphs, generalizing the main result of \cite{EKLP}.
\par James Propp, in \cite{Proppfaces}, proposed investigating the related {\em cube recurrence} given by $f_{i,j,k} = x_{i,j,k}$ for $i+j+k = -1, 0, 1$, and $$f_{i,j,k}f_{i-1,j-1,k-1} = f_{i-1,j,k}f_{i,j-1,k-1} + f_{i,j-1,k}f_{i-1,j,k-1} + f_{i,j,k-1}f_{i-1,j-1,k}.\quad (i + j + k > 1)$$  Fomin and Zelevinsky showed that the cube recurrence also generates Laurent polynomials.  Propp noticed empirically that each coefficient in these polynomials is equal to $1$ and that each variable takes only exponents in the range $-1,\ldots,4$.  Our goal is to construct combinatorial objects that are enumerated by a generalized form of the cube recurrence; Propp's observations, among other interesting results, will then follow directly.  Although the structure of the proof parallels that of Speyer for the octahedron recurrence, the combinatorial objects which we produce are quite different.
\par It is worth noticing that the cube recurrence can be written in a slightly more symmetric fashion: the families of functions $(f_{i,j,k})$ satisfying the cube recurrence can be made to correspond to the families $(g_{i,j,k})$ satisfying the recurrence $$g_{i,j,k}g_{i-1,j-1,k-1} + g_{i-1,j,k}g_{i,j-1,k-1} + g_{i,j-1,k}g_{i-1,j,k-1} + g_{i,j,k-1}g_{i-1,j-1,k} = 0$$ by taking $g_{i,j,k} = -f_{i,j,k}$ when $i+j+k \equiv 0 \mod{4}$, $g_{i,j,k} = f_{i,j,k}$ otherwise.  This latter equation has the esthetic advantage of being invariant not only under translation and permutation of coordinates but also under reflections (e.g. the substitution $i \leftarrow -i$).  However, we will not make further use of it here.

\section{The recurrence}
\label{recursec}
We will consider the polynomials generated by the cube recurrence using various sets of initial conditions.  In order to describe these initial conditions, we will need to develop some notation.
\par Define the {\em lower cone} of any $(i,j,k) \in \Z^3$ to be $$C(i,j,k) = \{(i',j',k') \in \Z^3\ |\ i' \leq i, j' \leq j, k' \leq k\}.$$  Let $\mL \subseteq \Z^3$ be a subset such that, whenever $(i,j,k) \in \mL, C(i,j,k) \subseteq \mL$.  (Thus, $\mL$ is an order-ideal in $\Z^3$ under the standard product ordering.)  Let $\mU = \Z^3 - \mL$, and define the set of {\em initial conditions} $$\mI = \{(i,j,k) \in \mL\ |\ (i+1,j+1,k+1) \in \mU\}.$$ To each $(i,j,k) \in \mI$ we assign a formal variable $x_{i,j,k}$.  We also define {\em edge variables} $a_{j,k}, b_{i,k}, c_{i,j}$ for all $i, j, k \in \Z$; the reason for this terminology will become clear later.
\par Now let $f_{i,j,k} = x_{i,j,k}$ for $(i,j,k) \in \mI$.  When $(i,j,k) \in \mU$ and $C(i,j,k) \cap \mU$ is finite, we define \begin{equation} \label{recur} f_{i,j,k} = \frac{b_{i,k}c_{i,j}f_{i-1,j,k}f_{i,j-1,k-1} + c_{i,j}a_{j,k}f_{i,j-1,k}f_{i-1,j,k-1} + a_{j,k}b_{i,k}f_{i,j,k-1}f_{i-1,j-1,k}}{f_{i-1,j-1,k-1}}. \end{equation}  (We leave $f_{i,j,k}$ undefined for all other points $(i,j,k)$.)  This recurrence gives us a well-defined rational function in the variables $\{a_{j',k'},$ $b_{i',k'}, c_{i',j'}, x_{i',j',k'}\ |\ (i',j',k') \in C(i,j,k) \cap \mI\}$, which takes a positive value when all the variables are set to $1$.  To see this, use induction on $|C(i,j,k) \cap \mU|$: either $(i-1,j,k) \in \mI$, or $(i-1,j,k) \in \mU$ and $|C(i-1,j,k) \cap \mU| < |C(i,j,k) \cap \mU|$; similarly for $(i,j-1,k-1)$, and so forth.  (We need the positivity statement in the induction hypothesis to ensure that the recurrence never produces a division by $0$.)  Also, to see that all edge variables appearing in $f_{i,j,k}$ really are of the form $a_{j',k'}, b_{i',k'},$ or $c_{j',k'}$ for some $(i',j',k') \in C(i,j,k) \cap \mI$, just notice (for example) that $(i,j,k') \in \mI$ for some $k' \leq k$: simply choose the maximal $k'$ for which $(i,j,k') \in \mL$, as, by finiteness, some such $k'$ must exist.
\par Henceforth, we will only investigate the value of $f_{i,j,k}$ for one particular $(i,j,k)$.  Because the definitions of $\mL, \mU, \mI$ and the recurrence itself are invariant under translation (modulo some relabeling of variables), we may assume that $(i,j,k) = (0,0,0)$.  We may also make another simplifying assumption:  Let $\mL' = \mL \cap C(0,0,0)$, and define $\mU', \mI'$ analogously to $\mU, \mI$.  It is easy to check that $\mI \cap C(0,0,0) \subseteq \mI'$.  In particular, this means that running the cube recurrence gives the same value for $f_{0,0,0}$ regardless of whether we use $\mI$ or $\mI'$ as our initial conditions, so we can safely replace $\mL$ by $\mL'$.  Therefore, we assume henceforward that $\mL \subseteq C(0,0,0)$ and that $C(0,0,0) \cap \mU$ is finite, except where explicitly stated otherwise.  
\par
\begin{figure}[hbt] \vskip 0.4in \begin{center} \includegraphics[width=5.5in]{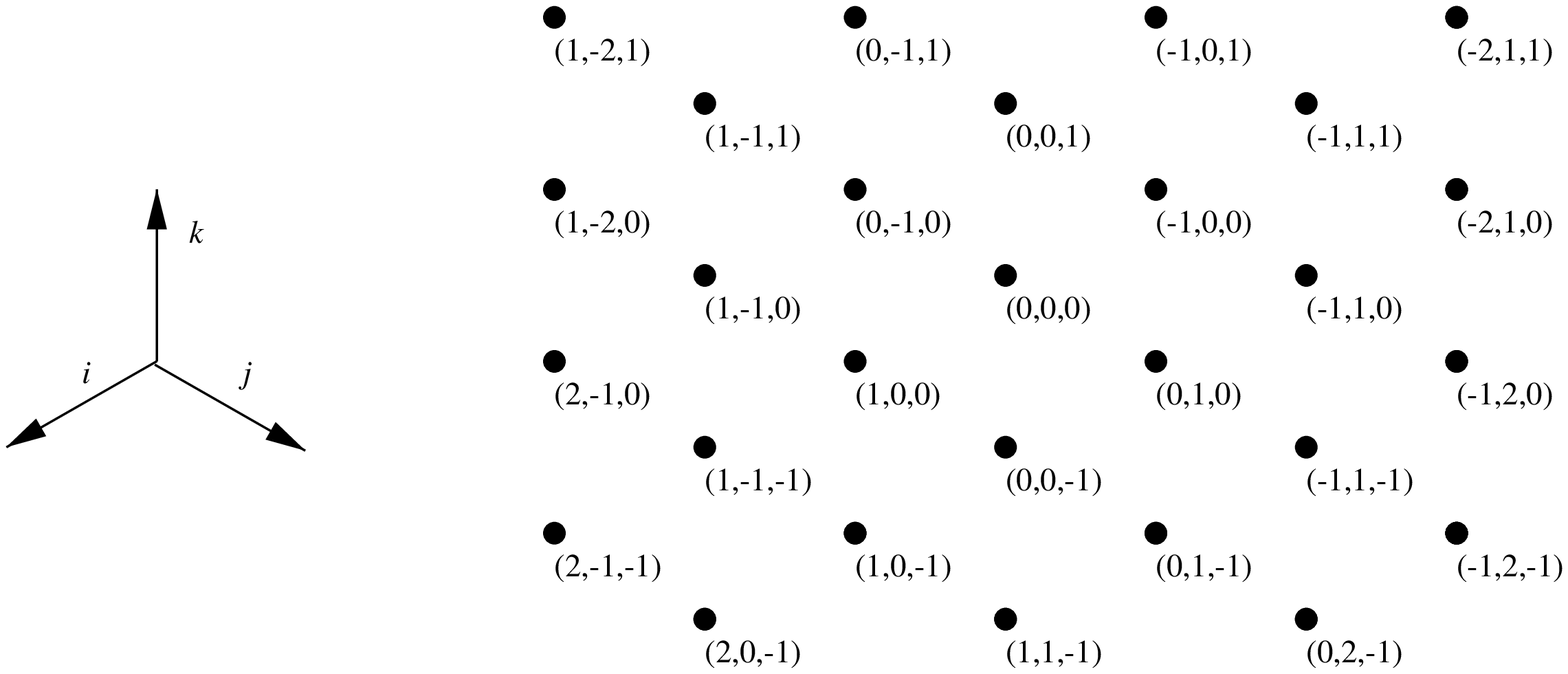} \caption{The triangular lattice} \label{tricoord}\end{center} \end{figure}

 We make one observation now that is not mathematically essential but will prove  extremely useful for visual intuition; it is based on the method of describing initial conditions used in \cite{Speyer}.  Let $S = \{(i,j,k) \in \Z^3\ |\ -1 \leq i+j+k \leq 1\}$.  For each $(i,j,k) \in S$, consider the largest integer $h$ such that $(i+h, j+h, k+h) \in \mL$.  (To see that such a value exists, notice that $(i+h,j+h,k+h) \in \mU$ for $h$ sufficiently large because $\mL \subseteq C(0,0,0)$; on the other hand, the finiteness condition ensures $(i+h,j+h,k+h) \in \mL$ for $h$ sufficiently negative.)  Then $(i+h',j+h',k+h') \in \mL$ for $h' \leq h$, while $(i+h',j+h',k+h') \in \mU$ for $h' > h$.  It follows that $h$ is the unique value such that $(i+h,j+h,k+h) \in \mI$.  The existence and uniqueness of such an $h$ for each $(i,j,k)$ imply that projecting $\mI$ onto the isometric plane $i+j+k = 0$ yields the same image as projecting $S$ onto this plane, namely the triangular lattice of Figure \ref{tricoord}.  Each point of the lattice is the image of exactly one point of $\mI$.
\par One further definition related to initial conditions will prove helpful.  By the finiteness assumption, we know that if $N$ is sufficiently large, we will have $(i,j,k) \in \mL$ whenever $i + j + k \leq -N$ (and $i,j,k \leq 0$).  In particular, we have that \begin{itemize} \item[(i)] $(i,j,k) \in \mI$ whenever $i + j + k \leq -N$ and $\max\{i,j,k\} = 0$; \item[(ii)] $(i,j,k) \notin \mI$ when $i + j + k \leq -N-3$ and $\max\{i,j,k\} < 0$. \end{itemize} Any $N$ with properties (i) and (ii) will be called a {\em cutoff} for $\mI$.

\section{Groves}
We will now introduce the combinatorial objects that will form the basis for our understanding of the cube recurrence.  We first provide the definition that will be most useful for purposes of subsequent proofs; later, we will offer an alternative representation that may be more practical as a kind of ``shorthand.''  We assume throughout that $\mL, \mU, \mI$ are as described in Section \ref{recursec}.
\par One preliminary notion that will prove crucial is that of a {\em rhombus}.  For each point $(i,j,k) \in \mI$, we define the following three sets: \begin{eqnarray*} r_a(i,j,k) &=& \{(i,j,k), (i,j-1,k), (i,j,k-1), (i,j-1,k-1)\} \\ r_b(i,j,k) &=& \{(i,j,k), (i-1,j,k), (i,j,k-1), (i-1,j,k-1)\} \\ r_c(i,j,k) &=& \{(i,j,k), (i-1,j,k), (i,j-1,k), (i-1,j-1,k)\} \end{eqnarray*}  We then define a rhombus to be any set of the form $r_a(i,j,k), r_b(i,j,k)$, or $r_c(i,j,k)$ that is contained in $\mI$.  Each rhombus can be decomposed into two pairs of points as follows: \begin{eqnarray*} e_a(i,j,k) = \{(i,j-1,k), (i,j,k-1)\} &\qquad& e'_a(i,j,k) = \{(i,j,k), (i,j-1,k-1)\} \\ e_b(i,j,k) = \{(i-1,j,k), (i,j,k-1)\} &\qquad& e'_b(i,j,k) = \{(i,j,k), (i-1,j,k-1)\} \\ e_c(i,j,k) = \{(i-1,j,k), (i,j-1,k)\} &\qquad& e'_c(i,j,k) = \{(i,j,k), (i-1,j-1,k)\} \end{eqnarray*} We refer to $e_a(i,j,k)$ as the {\em long edge} and $e'_a(i,j,k)$ as the {\em short edge} associated with $r_a(i,j,k)$, and similarly for the other pairs.
\par
\begin{figure}[hbt] \begin{center} \includegraphics[width=3in]{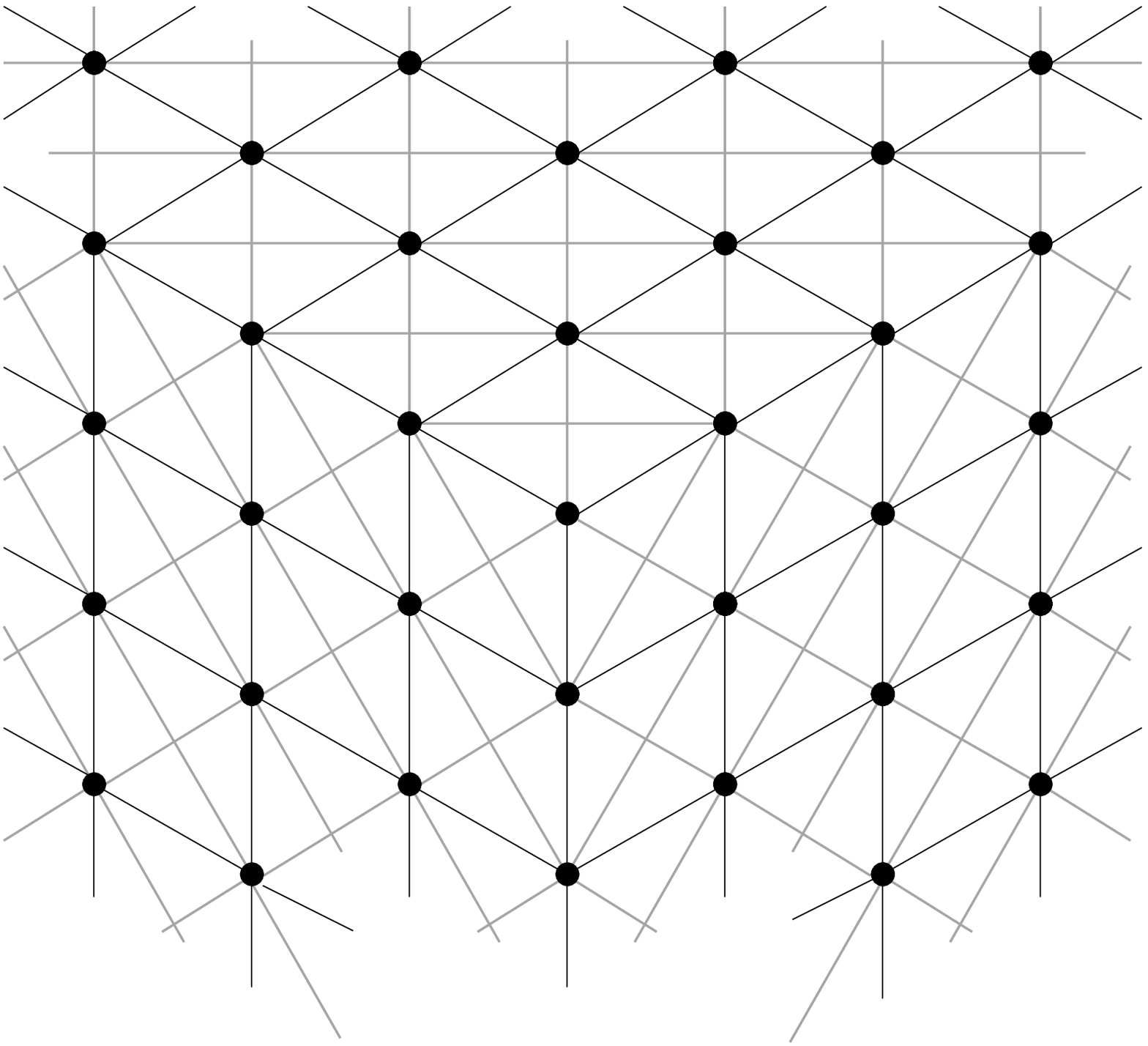} \caption{$\mG$ for the initial conditions induced by $\mL = C(0,0,0)$} \label{trivmg} \end{center} \end{figure}

 We now construct an (infinite) graph $\mG$ whose vertices are the points in $\mI$.  The edges of $\mG$ are simply the long and short edges of all rhombi occurring in $\mI$.  The motivation for our terminology becomes clear when we project $\mI$ onto the plane $i+j+k = 0$.  An example is shown in Figure \ref{trivmg}, where $\mL$ is simply all of $C(0,0,0)$; the rhombi are outlined in thin solid black, and the gray lines represent edges of $\mG$.
\par In fact, any choice of $\mI$ produces a rhombus tiling of the plane; we do not prove this here, since it is ancillary to our main concerns, although it will become apparent from techniques to be introduced subsequently.
\par
\begin{figure}[hbt] \begin{center} \includegraphics[width=4in]{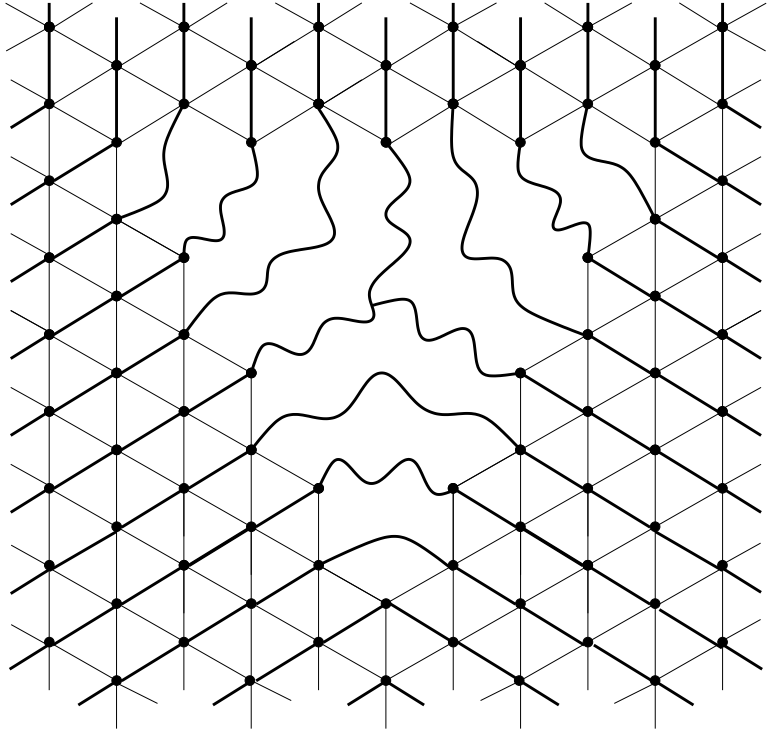} \caption{Compactness and connectivity conditions for a generic grove, $N = 3$} \label{conncond} \end{center} \end{figure}

\begin{figure}[htb] \begin{center} \includegraphics[width=4in]{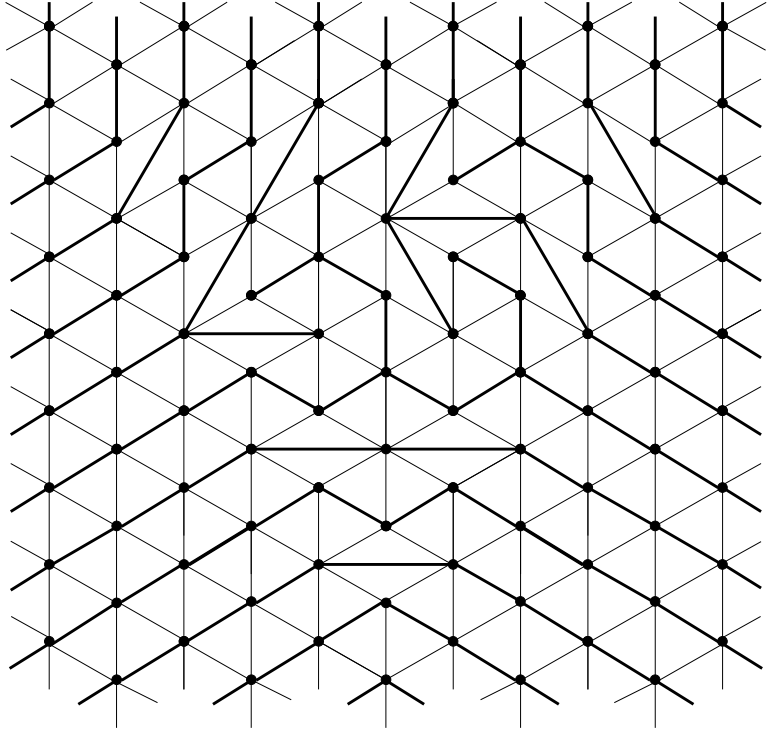} \caption{An example of a grove} \label{stdgrov} \end{center} \end{figure}

 Now suppose that $N$ is a cutoff for $\mI$. We define an {\em $\mI$-grove within radius $N$} to be a subgraph $G \subseteq \mG$ with the following properties:
\begin{itemize}
\item (Completeness) the vertex set of $G$ is all of $\mI$;
\item (Complementarity) for every rhombus, exactly one of its two edges occurs in $G$;
\item (Compactness) for every rhombus all of whose vertices satisfy $i + j + k < -N$, the short edge occurs in $G$;
\item (Connectivity) every component of $G$ contains exactly one of the following sets of vertices, and conversely, each such set is contained in some component: \begin{itemize} \item $\{(0, p, q), (p, 0, q)\}, \{(p, q, 0), (0, q, p)\}$, and $\{(q, 0, p), (q, p, 0)\}$ for all $p, q$ with $0 > p > q$ and $p + q \in \{-N-1, -N-2\}$; \item $\{(0,p,p), (p,0,p), (p,p,0)\}$ for $2p \in \{-N-1, -N-2\}$; \item $\{(0,0,q)\}, \{(0,q,0)\}$, and $\{(q,0,0)\}$ for $q \leq -N-1$. \end{itemize}
\end{itemize}
Loosely speaking, then, a grove consists of a fixed set of edges outside the region $\{(i,j,k) \in \mI\ |\ i+j+k \geq -N\}$, together with a graph inside this region constrained by connectivity conditions among the vertices near the boundary of the region.  Figure \ref{conncond} illustrates this structure (here for $N = 3$).  The thin lines again signify rhombus boundaries; the thick black lines are the short edges forced by compactness, and the wavy lines connect vertices that must belong to the same component in any grove.  Figure \ref{stdgrov} shows an example of an actual grove within radius $4$ on the initial conditions obtained by taking $\mL = \{(i,j,k) \in C(0,0,0)\ |\ i + j + k \leq -4\}$; the thick black lines are just the edges of $G$.
\par We will show that the definition of a grove actually does not depend on the choice of $N$ (as long as $N$ is a cutoff for $\mI$).  Although it is possible to prove this directly from the definition, we do not do so here, as the verifications are somewhat tedious; instead, we will obtain it as a consequence of our main theorem.
\par In order to state the theorem, it will be necessary to indicate how groves can be represented algebraically.  Given a grove $G$ within radius $N$, define the corresponding Laurent monomial $$m(G) = \left( \product_{e_a(i,j,k) \in E(G)} a_{j,k} \right) \left( \product_{e_b(i,j,k) \in E(G)} b_{i,k} \right) \left( \product_{e_c(i,j,k) \in E(G)} c_{i,j} \right) \left( \product_{(i,j,k) \in \mI} x_{i,j,k}^{\deg(i,j,k) - 2} \right).$$  The first three products are finite because they are simply products of edge variables corresponding to long edges in $G$, and the compactness condition ensures that only finitely many long edges appear.  (These products also elucidate our use of the term ``edge variable.'')  The last product is finite because the compactness condition ensures that each $(i,j,k) \in \mI$ with $i+j+k \leq -N-3$ has degree $2$: since $N$ is a cutoff, $\max\{i,j,k\} = 0$; without loss of generality, let $i \leq j \leq k = 0$, and then either $j = 0$, and $(i,0,0)$ is adjacent to $(i-1,-1,0)$ and $(i-1,0,-1)$ (and no other vertices), or $j < 0$, and $(i,j,0)$ is adjacent to $(i+1,j+1,0)$ and $(i-1,j-1,0)$.  (See e.g. figure \ref{conncond}.)  There are then only finitely many $(i,j,k)$ with $i+j+k > -N-3$, and only the $x_{i,j,k}$ corresponding to these points can contribute to the fourth product in $m(G)$.  
\par Notice also that $m(G)$ uniquely determines $G$ (independently of $N$), since it states precisely which long edges occur in $G$.  To see this, we need only observe that no two distinct long edges can be represented by the same edge variable.  For example, if $a_{j,k}$ represented both $e_a(i,j,k)$ and $e_a(i',j,k)$ with $i > i'$, then $(i,j,k) \in r_a(i,j,k) \subseteq \mI$, so $(i'+1,j,k) \in \mL$ and $(i',j-1,k-1) \notin \mI$, but also $(i',j-1,k-1) \in r_a(i',j,k) \subseteq \mI$, a contradiction.
\par At this point we are prepared to state our main theorem.
\begin{theorem} \label{mainthm} Let $\mL, \mU, \mI$ be as described in Section \ref{recursec}, including the assumptions $\mL \subseteq C(0,0,0)$ and $\mU \cap C(0,0,0)$ finite.  Define $f_{i,j,k}$ as in Section \ref{recursec}, and let $N$ be any cutoff for $\mI$.  Then $$f_{0,0,0} = \sum_G m(G),$$ where the sum is taken over all $\mI$-groves within radius $N$. \end{theorem}
\par The proof is postponed for the purpose of generating suspense.  We point out now, however, that there is only one way of decomposing $f_{0,0,0}$ as a sum of Laurent monomials in the variables $\{a_{j,k},b_{i,k},c_{i,j},x_{i,j,k}\}$.  Because each monomial $m(G)$ in turn determines $G$ uniquely, we see that we obtain the same set of radius-$N$ groves regardless of the choice of cutoff $N$.  Consequently, in all subsequent discussion (except the proof of Theorem \ref{mainthm} itself), we may drop the $N$ and simply use the term ``grove,'' or ``$\mI$-grove'' if the choice of $\mI$ is ambiguous.
\par Before proceeding, we make one basic observation about the structure of groves that is not necessarily apparent from the definition. \begin{theorem} \label{acyclic} Every grove is acyclic. \end{theorem}  The proof will require one preliminary result.  Let $\mJ = \{(i,j,k) \in \mI\ |\ i + j + k \geq -N-2\}$.  We then have \begin{lemma} \label{rhombcount} The set $\mJ$ consists of $3{N+3\choose 2} + 1$ points and contains exactly $3{N+2\choose 2}$ rhombi. \end{lemma}  The proof is deferred, as it will use the same technique as the proof of Theorem \ref{mainthm}.
\pfof{Theorem \ref{acyclic}} Let $G$ be a grove on the initial conditions $\mI$; we begin by proving that $H$, the induced subgraph on $\mJ$, is acyclic.  We claim that if any two vertices of $\mJ$ lie in the same component of $G$, they are connected by a path contained in $\mJ$.  For suppose not; choose two vertices connected by a path (which we may assume to have no repeated vertices) not contained in $\mJ$.  Choose a vertex $(i,j,k)$ of this path with $i + j + k$ minimal; in particular, $i + j + k < -N - 2$, and $(i,j,k)$ is not an endpoint of the path.  Because $N$ is a cutoff, $\max\{i,j,k\} = 0$; assume without loss of generality $i = 0$.  Then $(0,j,k)$ belongs only to the rhombi $r_a(0,j+1,k), r_a(0,j,k+1)$, each of which (by compactness) contributes to $G$ the edge not incident to $(0,j,k)$; $r_a(0,j,k)$, which contributes the edge $e'_a(0,j,k) = \{(0,j,k), (0,j-1,k-1)\}$, which cannot be used in the path (by minimality); and $r_a(0,j+1,k+1)$, which contributes the edge $e'_a(0,j+1,k+1) = \{(0,j+1,k+1), (0,j,k)\}$.  Hence, only one edge incident to $(0,j,k)$ may be used in the path, contradicting the assumption of no repeated vertices.  The claim follows.  Hence, every component of $G$ that contains any vertex of $\mJ$ induces a single component of $H$.
\par Now consider the following classes of vertices in $\mJ$: \begin{itemize} \item $\{(0, p, q), (p, 0, q)\}, \{(p, q, 0), (0, q, p)\}$, and $\{(q, 0, p), (q, p, 0)\}$ for all $p, q$ with $0 > p > q$ and $p + q \in \{-N-1, -N-2\}$; \item $\{(0,p,p), (p,0,p), (p,p,0)\}$ for $2p \in \{-N-1, -N-2\}$; \item $\{(0,0,q)\}, \{(0,q,0)\}$, and $\{(q,0,0)\}$ for $q \in \{-N-1, -N-2\}$. \end{itemize}  By the foregoing and connectivity for $G$, each class is contained in a single component of $H$, and certainly no component may contain vertices of more than one class (otherwise the corresponding component of $G$ would, which is impossible).  We also claim that every component of $H$ contains one of the above classes of vertices; it suffices to show that the corresponding component of $G$ contains some vertex $(i,j,k)$ with $i + j + k \in \{-N-1, -N-2\}$.  If not, then by the connectivity condition, this component of $G$ must contain some vertex $(i,j,k)$ with $i + j + k < -N-2$.  But it is impossible for such a vertex to be connected to a vertex of $\mJ$ by a path not going through any point $(i,j,k)$ with $i + j + k \in \{-N-1, -N-2\}$, since the sum of coordinates changes by at most two at each step along the path.  This is a contradiction.
\par We therefore conclude that the components of $H$ are in bijection with our classes of vertices, of which there are $3N + 7$.  On the other hand, by Lemma \ref{rhombcount} (and the complementarity condition), $H$ has $3{N+3 \choose 2} + 1$ vertices and $3{N+2 \choose 2}$ edges, for a minimum of $$\left(3{N+3 \choose 2} + 1\right) - 3{N+2 \choose 2} = 3N + 7$$ components, with equality only if $H$ is acyclic.  Equality does hold, so $H$ is acyclic, as claimed.
\par Now, suppose $G$ contains some cycle.  Since the definition of a grove is independent of the choice of cutoff $N$, we may choose $N$ large enough so that all the vertices of the cycle belong to $\mJ$.  Then the induced subgraph $H$ contains a cycle, and this is a contradiction.  Hence, $G$ is acyclic. \eop
\par The definition of a grove we have presented is somewhat cumbersome.  For purposes of empirical investigation, infinite graphs are inconvenient to work with; groves as defined above also contain, in a certain sense, redundant information.  We therefore will present --- in slightly less detail --- a simplification which we have at times found intuitively more useful, in the hope that it will also prove more practical for later investigations.
\par By projecting onto the plane $i + j + k = 0$, we can reinterpret groves as graphs on the infinite triangular lattice.  It is not hard to check that the interiors of distinct rhombi cannot overlap in this projection; this fact, in conjunction with the complementarity requirement, ensures that the resulting graphs are planar.  
\par We can define a point $(i,j,k) \in \Z^3$ to be {\em even} or {\em odd} depending on the parity of $i+j+k$.  Every rhombus then has two even vertices and two odd vertices, and the edge used in $G$ connects two vertices of the same parity (and so can also be called {\em even} or {\em odd}).  It follows that, given ``half'' of a grove, we can uniquely reconstruct the other half: if we know which even edges are used in $G$, the remaining edges must be precisely the odd edges of those rhombi whose even edges are not used.
\par
\begin{figure}[bthp] \begin{center} \includegraphics[width=6in]{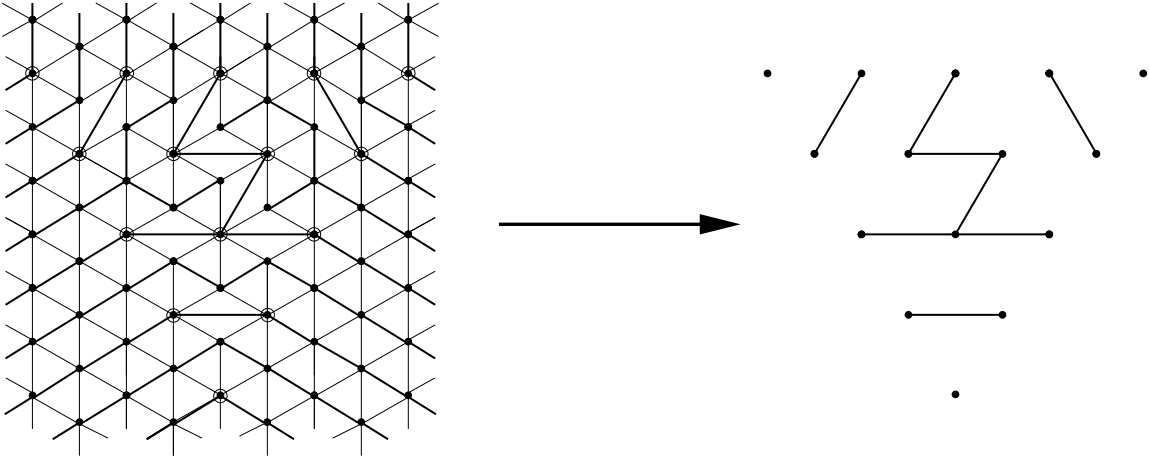} \caption{Obtaining a simplified grove from a grove} \end{center} \label{simpgrov} \end{figure}

 Finally, we can restrict our attention to a finite piece of $\mI$, since, whenever a vertex $(i,j,k)$ satisfies $i+j+k \leq -N-3$, we know precisely which edges are incident to it, by the compactness condition.  Putting all these observations together, we define a {\em simplified grove within radius $N$}, where $N$ is a cutoff for $\mI$ and furthermore is odd, to be a subgraph $G'$ of $\mG$ satisfying:
\begin{itemize}
\item (Vertex set) the vertex set of $G'$ is $\{(i,j,k) \in \mI\ |\ i + j + k \equiv 0 \mod{2};\ i+j+k \geq -N-1\}$;
\item (Acyclicity) $G'$ is acyclic;
\item (Connectivity) the {\em even boundary vertices} $\{(i,j,k) \in \mI\ |\ i + j + k = -N - 1;$ $\max\{i,j,k\} = 0\}$ can be partitioned into the following sets so that each component of $G'$ contains exactly one set, and conversely, each set is contained in some component: \begin{itemize} \item $\{(0,p,q), (p,0,q)\}, \{(p,q,0), (0,q,p)\}$, and $\{(q,p,0), (q,0,p)\}$ for $0 > p > q$, $p + q = -N - 1$; \item $\{(0,\frac{-N-1}{2},\frac{-N-1}{2}), (\frac{-N-1}{2}, 0, \frac{-N-1}{2}), (\frac{-N-1}{2}, \frac{-N-1}{2}, 0)\}$, \item $\{(0,0,-N-1)\}, \{(0,-N-1,0)\}$, and $\{(-N-1,0,0)\}$. \end{itemize}
\end{itemize}
For any grove $G$, the induced subgraph on $\{(i,j,k)\ |\ i+j+k \equiv 0 \mod 2;\ i+j+k \geq -N-1\}$ is then a simplified grove.  (An example of a grove and the corresponding simplified grove, with $N = 3$, is shown in Figure \ref{simpgrov}; the vertices circled on the left side of the figure are those belonging to the vertex set of the simplified grove.)  The only condition that is not entirely evident is connectivity, which follows from the connectivity condition on $G$, with a slight bit of subtlety: we must check that, in the grove $G$, a path connecting any two even boundary vertices as specified by connectivity is in fact contained in the vertex set of $G'$.  However, the requirement that all vertices on the path be even poses no difficulty, as even vertices are connected only to even vertices; and the verification that all vertices $(i,j,k)$ on the path satisfy $i+j+k \geq -N-1$ is precisely as described in the proof of Theorem \ref{acyclic}.
\par We next show the converse of the above --- every simplified grove is induced by a unique grove.  Given a simplified grove $G'$, we again let $\mJ = \{(i,j,k) \in \mI\ |\ i + j + k \geq -N-2\}$.  We extend $G'$ to a graph $G$ on all of $\mI$ as follows: for every rhombus contained in $\mJ$ whose even edge does not appear in $G'$, we include the odd edge instead; we then include the short edges of all other rhombi in $\mI$.  In view of the complementarity and compactness conditions, the resulting graph is the only possible grove on $\mI$ that induces the simplified grove $G'$.  We claim that this $G$ is indeed a grove.  All the conditions except connectivity are immediate; this last is a bit more involved.
\par We first show that our graph $H$ on the vertex set $\mJ$, obtained by including the appropriate odd edges, is acyclic.  It cannot have any cycles on even vertices, as such a cycle would already have existed in $G'$.  If there is a cycle on the odd vertices, then it is not hard to see that its planar projection must enclose some even vertices of $\mJ$ but cannot enclose any even boundary vertices.  By planarity, then, $G'$ has some vertices that cannot be connected to any even boundary vertices, violating its connectivity requirement.  Thus, no cycles have been introduced in $\mJ$.
\par Now we again apply the component-counting technique.  From connectivity for a simplified grove, we can check that there are always $(3N+5)/2$ components.  If we consider the graph $H$, Lemma \ref{rhombcount} tells us that we have $3{N+2\choose 2}$ edges and $3{N+3\choose 2}+1$ vertices; since the graph is acyclic, it has $$3{N+3\choose 2}+1 - 3{N+2\choose 2} = 3N+7$$ components.  In particular, there are $(3N+9)/2$ components on the odd vertices.
\par Now divide the {\em odd boundary vertices} $\{(i,j,k)\ |\ i+j+k = -N - 2;\ \max\{i,j,k\} = 0\}$ into classes according to connectivity for a grove: we have the classes $\{(0,p,q), (p,0,q)\}$, $\{(p,q,0), (0,q,p)\}$, $\{(q,p,0), (q,0,p)\}$ for $0 > p > q, p + q = -N-2$, as well as $\{(0,0,-N-2)\}, \{(0,-N-2,0)\}, \{(-N-2,0,0)\}$.  We thus get $(3N+9)/2$ classes.  Moreover, we claim that, in our graph on $\mI$, no two odd boundary vertices from different classes can lie in the same component.  To see this, we consider the paths of $G'$ connecting the even boundary vertices, together with the extra short edges emanating from them (given by the compactness condition), and project into the plane.  These paths divide the plane into sectors, each of which contains the odd boundary vertices of (at most) one class.  If any two odd boundary vertices from different classes lie in the same component of $G$, the path connecting them must intersect one of the paths on even vertices.  However, the two paths cannot intersect at an actual vertex, since one path uses only even vertices and the other uses only odd vertices.  This violates planarity.
\par
\begin{figure}[htb] \begin{center} \includegraphics[width=4in]{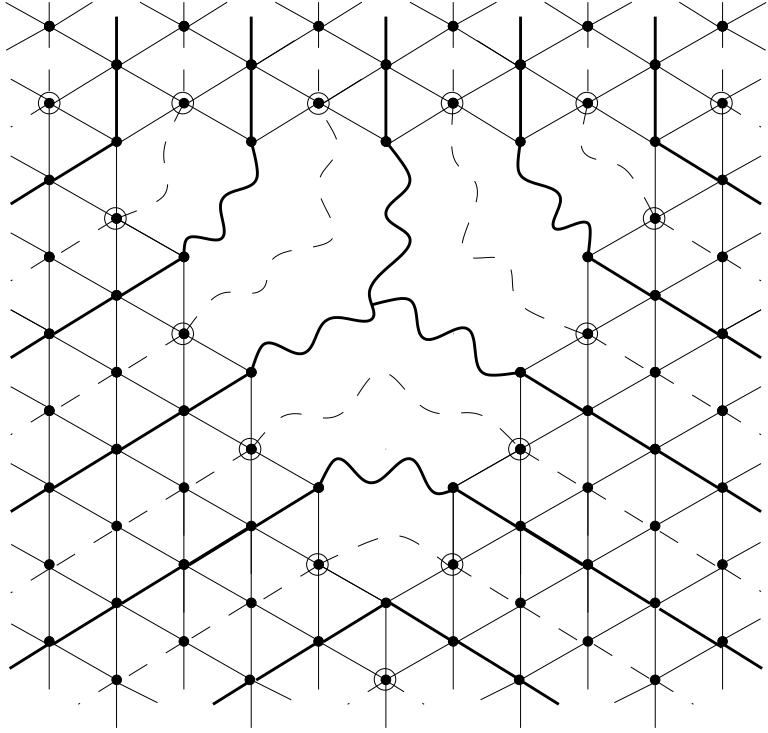} \caption{Connectivity restrictions among odd boundary vertices} \label{esectors} \end{center} \end{figure}

 This argument is shown in Figure \ref{esectors}.  The solid wavy lines represent paths between even boundary vertices, as specified by simplified-grove connectivity; the odd boundary vertices are circled, and the dashed wavy lines show all possible connectivities that can occur between them without violating planarity.  These are precisely the connectivities required for a grove. 
\par At this point we have shown that, for any component of $G$ --- and so for any component of $H$ --- the odd boundary vertices it contains all belong to the same class.  Since we have $(3N+9)/2$ classes, we get at least $(3N+9)/2$ components of $H$ on odd vertices.  But we have already shown that equality occurs.  This is only possible if, for each class, all its vertices belong to the same component, and every component contains one such class of vertices.  We may now conclude that each class of vertices required for a grove is contained in some component of $G$: either the class consists of odd boundary vertices, and we have just proven that they are connected; it consists of even boundary vertices, and the conclusion follows from connectivity for $G'$; or it consists of a single vertex outside of $\mJ$, and the conclusion is trivial.  Conversely, consider any component of $G$.  If it contains any vertex of $\mJ$, it either contains odd vertices, and then (as we have just seen) it contains an entire class of odd boundary vertices, or it contains an even vertex, and then (by connectivity for $G'$) it contains an entire class of even boundary vertices.  Moreover, no two boundary vertices from different classes may be connected within $\mJ$, and the extra short edges cannot introduce any new such connectivities (by the argument in the proof of Theorem \ref{acyclic}); thus, no component of $G$ contains more than one class of boundary vertices.  On the other hand, if our component of $G$ does not contain any vertex of $\mJ$, then it is a chain of short edges and contains exactly one vertex of the form $(0,0,q)$, $(0,q,0)$, or $(q,0,0)$ ($q \leq -N-3$).  Every such vertex lies in a component of $G$ consisting entirely of short edges, which cannot intersect $\mJ$ and therefore contains none of the classes of boundary vertices.  Thus, each component of $G$ contains exactly one class of vertices.  We now have verified that the grove connectivity condition is met in its entirety.
\par This completes the proof that every simplified grove is induced by a unique grove.  We now have established a bijection between groves and simplified groves (within radius $N$), so the latter may be used as a sort of shorthand to represent the former.

\section{Some consequences}
Before proceeding to the proof of the main theorem, let us continue to assume that it holds and note some consequences.  As usual, the hypotheses of Section \ref{recursec} hold, including the assumptions that $\mL \subseteq C(0,0,0)$ and $|\mU \cap C(0,0,0)| < \infty$, except where stated otherwise.
\par It is clear that, as a Laurent polynomial in $\{a_{j,k}, b_{i,k}, c_{i,j}, x_{i,j,k}\}$, $f_{0,0,0}$ has every coefficient equal to $1$ --- that is, a grove $G$ is uniquely determined by the monomial $m(G)$, as already discussed.  In \cite{Proppfaces}, Propp conjectured that, for the initial conditions $\mI = \{(i,j,k) \in \Z^3\ |\ -1 \leq i+j+k \leq 1\}$, every coefficient is still $1$ if we set all edge variables equal to 1 and simply view each $f_{i,j,k}$ as a Laurent polynomial in the variables $x_{i',j',k'}$.  Fomin and Zelevinsky (in \cite{FomZel}) conjectured that the coefficients remain at least nonnegative for general initial conditions.  We will show that the coefficients are all $1$ when all edge variables are set to $1$ and the initial conditions are arbitrary.  This is equivalent to the assertion that any grove $G$ (on a given set of initial conditions) can be uniquely reconstructed from the degrees of its vertices.
\par More precisely, we will prove:
\begin{theorem} \label{coeffone} Suppose $r_a(i_0,j_0,k_0) \subseteq \mI$ is a rhombus.  Then, for any grove $G$, $$\sum_{^{(i,j,k) \in \mI}_{j < j_0, k < k_0}} \left( \deg(i,j,k) - 2 \right)$$ equals $-1$ if the long edge $e_a(i_0,j_0,k_0)$ appears in $G$, and $0$ if it does not appear. \end{theorem}
(We have already noted that the sum has only finitely many nonzero terms.)  This theorem, and the analogous statements for the other two types of rhombi, will imply Propp's conjecture for any $\mI$.
\pf Let $t$ be a new formal indeterminate.  We apply the following variable substitutions: $$a_{j_0,k_0} \leftarrow t; \qquad a_{j,k} \leftarrow 1 \hbox{ otherwise};$$ $$b_{i,k}, c_{i,j} \leftarrow 1;$$ $$x_{i,j,k} \leftarrow t\ (j < j_0, k < k_0); \qquad x_{i,j,k} \leftarrow 1 \hbox{ otherwise}.$$  We proceed to compute values of $f_{i,j,k}$ $(i,j,k \in \mI \cup \mU)$ by the recurrence (\ref{recur}); each such value will then be a Laurent polynomial in $t$.  We claim that, in fact, $f_{i,j,k}$ is a constant multiple of $t$ for $j < j_0, k < k_0$, and is simply a constant otherwise.  (Notice also that these constants are always positive.)
\par The proof is by induction on the cardinality of $C(i,j,k) \cap \mU$.  If this cardinality is $0$, we have $(i,j,k) \in \mI$, and $f_{i,j,k} = x_{i,j,k} = t$ or $1$ according to whether or not $j < j_0$ and $k < k_0$.  Otherwise, we split into cases:
\begin{itemize} \item If $j < j_0, k < k_0$, then (\ref{recur}) simply says that $$f_{i,j,k} = (f_{i-1,j,k}f_{i,j-1,k-1} + f_{i,j-1,k}f_{i-1,j,k-1} + f_{i,j,k-1}f_{i-1,j-1,k}) / f_{i-1,j-1,k-1}.$$  However, either $(i-1,j,k) \in \mI$, or $(i-1,j,k) \in \mU$ and $|C(i-1,j,k) \cap \mU| < |C(i,j,k) \cap \mU|$ (this was discussed in Section \ref{recursec}), so the induction hypothesis tells us that $f_{i-1,j,k}$ is a constant multiple of $t$.  Similarly, $f_{i,j-1,k-1}, f_{i,j-1,k}$, and so forth are all constant multiples of $t$, and we conclude that $f_{i,j,k}$ is as well.
\item If $j = j_0, k < k_0$, then again $$f_{i,j,k} = (f_{i-1,j,k}f_{i,j-1,k-1} + f_{i,j-1,k}f_{i-1,j,k-1} + f_{i,j,k-1}f_{i-1,j-1,k}) / f_{i-1,j-1,k-1}.$$  The induction hypothesis now tells us that $f_{i-1,j,k}, f_{i-1,j,k-1}, f_{i,j,k-1}$ are constants, while $f_{i,j-1,k-1}, f_{i,j-1,k}, f_{i-1,j-1,k}, f_{i-1,j-1,k-1}$ are constant multiples of $t$, and it follows that $f_{i,j,k}$ is also a constant.
\item If $j < j_0, k = k_0$, the reasoning is the same as in the previous case.
\item If $j = j_0, k = k_0$, then (\ref{recur}) takes the form $$f_{i,j,k} = (f_{i-1,j,k}f_{i,j-1,k-1} + tf_{i,j-1,k}f_{i-1,j,k-1} + tf_{i,j,k-1}f_{i-1,j-1,k}) / f_{i-1,j-1,k-1}.$$  The induction hypothesis now tells us that $f_{i-1,j,k}, f_{i,j-1,k}, f_{i,j,k-1}, f_{i-1,j,k-1}, f_{i-1,j-1,k}$ are constants, while $f_{i,j-1,k-1}, f_{i-1,j-1,k-1}$ are constant multiples of $t$.  We conclude that $f_{i,j,k}$ is again a constant.
\item Finally, if $j > j_0$ or $k > k_0$, then the recurrence again takes the form $$f_{i,j,k} = (f_{i-1,j,k}f_{i,j-1,k-1} + f_{i,j-1,k}f_{i-1,j,k-1} + f_{i,j,k-1}f_{i-1,j-1,k}) / f_{i-1,j-1,k-1}.$$  All the $f$'s appearing on the right are now constants, so $f_{i,j,k}$ is as well.
\end{itemize}
\par This completes the induction.  In particular, we know now that $f_{0,0,0}$ is a constant polynomial.  So, for every grove $G$, $m(G)$ becomes a constant under our variable substitutions --- the total exponent of $t$ is $0$.  This means that either the sum in the problem statement is $-1$ and the long edge $e_a(i_0,j_0,k_0)$ (corresponding to $a_{j_0,k_0}$) does appear, or the sum is $0$ and the long edge does not appear, as claimed. \eop
\par Another conjecture appearing in \cite{Proppfaces} is that, in each term of any polynomial generated by the cube recurrence, every $x_{i,j,k}$ has its exponent in the range $\{-1,0,\ldots,4\}$.  In terms of groves, this is equivalent to the assertion that every vertex has degree no less than $1$ and no more than $6$.  The lower bound is obvious, since the connectivity condition ensures that there are no isolated vertices (except possibly those of the forms $(0,0,q), (0,q,0), (q,0,0)$, but these lie on short edges given by the compactness condition).  The upper bound holds because each vertex has at most six neighbors in $\mG$; this essentially follows from the fact that, in a rhombus tiling of the plane, no vertex can belong to more than six rhombi.  To make the argument algebraically precise, we consider any vertex $(i,j,k) \in \mI$.  A priori, there are twelve vertices to which $(i,j,k)$ could be adjacent: those of the form $(i \pm 1, j \pm 1, k), (i \pm 1, j, k \pm 1)$, or $(i, j \pm 1, k \pm 1)$.  However, we will classify these vertices into six pairs such that, in each pair, only one of the two vertices can actually be adjacent to $(i,j,k)$ in $\mG$.  For example, if $(i,j,k)$ were adjacent to both $(i,j-1,k-1)$ and $(i+1,j,k-1)$, the latter adjacency would require the presence of the rhombus $r_b(i+1,j,k) \subseteq \mI$, but we cannot have $(i+1,j,k), (i,j-1,k-1)$ both in $\mI$.  The full pairing is shown in the Table \ref{sixrhombs}.
\begin{table}[hbt] \small \begin{tabular}{c|c|c}
For $(i,j,k)$ to have these neighbors & $\mI$ must contain this rhombus & causing this contradiction \\ \hline
$(i+1,j,k+1), (i+1,j-1,k)$ & $r_c(i+1,j,k)$ & $(i+1,j,k+1), (i,j-1,k)$ \\
$(i,j-1,k-1), (i+1,j,k-1)$ & $r_b(i+1,j,k)$ & $(i+1,j,k), (i,j-1,k-1)$ \\
$(i+1,j+1,k), (i,j+1,k-1)$ & $r_a(i,j+1,k)$ & $(i+1,j+1,k), (i,j,k-1)$ \\
$(i-1,j,k-1), (i-1,j+1,k)$ & $r_c(i,j+1,k)$ & $(i,j+1,k), (i-1,j,k-1)$ \\
$(i,j+1,k+1), (i-1,j,k+1)$ & $r_b(i,j,k+1)$ & $(i,j+1,k+1), (i-1,j,k)$ \\
$(i-1,j-1,k), (i,j-1,k+1)$ & $r_a(i,j,k+1)$ & $(i,j,k+1), (i-1,j-1,k)$
\end{tabular} \normalsize \label{sixrhombs} \caption{Proof that every vertex has degree $\leq 6$} \end{table}
\par This shows that every vertex has at most $6$ neighbors in $\mG$ (and so in $G$), as claimed.
\par The preceding results were originally phrased in \cite{Proppfaces} as algebraic statements about the cube recurrence, but they can be interpreted as geometric statements about the structure of groves.  Another geometric fact worth noting concerns the distribution of the orientations of long edges in any grove.
\begin{theorem} \label{triineq} Suppose that a grove $G$ has $n_a$ long edges of the form $e_a(i,j,k)$, $n_b$ of the form $e_b(i,j,k)$, and $n_c$ of the form $e_c(i,j,k)$.  Then the numbers $n_a + n_b - n_c, n_b + n_c - n_a, n_c + n_a - n_b$ are all nonnegative and even. \end{theorem}
\pf The proof is similar to that used for Theorem \ref{coeffone}.  We view $n_a, n_b, n_c$ as functions of $G$; it suffices to prove the result for $n_a(G) + n_b(G) - n_c(G)$, as the other cases are analogous.  We let $t$ be a formal indeterminate and make the variable substitutions \begin{eqnarray*} a_{j,k}, b_{i,k} &\leftarrow& t; \\ c_{i,j} &\leftarrow& 1/t; \\ x_{i,j,k} &\leftarrow& 1. \end{eqnarray*}  We again compute values of $f_{i,j,k}$ by recurrence (\ref{recur}); we know from Theorem \ref{mainthm} that each $f_{i,j,k}$ is a Laurent polynomial in $t$ with positive coefficients.  We claim that $f_{i,j,k}$ is in fact a polynomial in $t^2$ whose constant term is positive.
\par The proof is by induction on $|C(i,j,k) \cap \mU|$, as usual.  If this cardinality is $0$, then $(i,j,k) \in \mI$, so $f_{i,j,k} = x_{i,j,k} = 1$, and the result holds.  Otherwise, the recurrence (\ref{recur}) tells us that $$f_{i,j,k}f_{i-1,j-1,k-1} = f_{i-1,j,k}f_{i,j-1,k-1} + f_{i,j-1,k}f_{i-1,j,k-1} + t^2f_{i,j,k-1}f_{i-1,j-1,k}.$$  The induction hypothesis ensures that the right side is a polynomial in $t^2$ with positive constant term (since each factor is).  We also know that $f_{i-1,j-1,k-1}$ is a polynomial in $t^2$ with positive constant term.  Let $d$ be the minimal exponent for which the coefficient of $t^d$ in $f_{i,j,k}$ is nonzero; then the lowest power of $t$ appearing on the right side must also be $t^d$, so $d = 0$.  This assures us that $f_{i,j,k}$ is a polynomial in $t$ with positive constant term.  However, since $f_{i-1,j,k}$ and so forth are actually polynomials in $t^2$, $f_{i,j,k}$ is a rational function of $t^2$; hence, it is a polynomial in $t^2$, and the induction step holds.
\par In particular, $\sum_G m(G) = f_{0,0,0}$ is a polynomial in $t^2$.  On the other hand, for each grove $G$, it is apparent that $m(G) = t^{n_a(G) + n_b(G) - n_c(G)}$.  It follows that $n_a(G) + n_b(G) - n_c(G)$ is nonnegative and even. \eop
\par This observation allows us to provide a complete combinatorial proof of the Laurent property of the cube recurrence as stated (and proved using cluster-algebra techniques) by Fomin and Zelevinsky in \cite{FomZel}.  Their version is essentially as follows:
\begin{theorem} \label{fzcube} Let $\mL \subseteq \Z^3$ such that whenever $(i,j,k) \in \mL$ and $i' \leq i, j' \leq j, k' \leq k$, we also have $(i',j',k') \in \mL$.  Let $\mU = \Z^3 - \mL, \mI = \{(i,j,k) \in \mL\ |\ (i+1,j+1,k+1) \in \mU\}$, and suppose $C(i,j,k) \cap \mU$ is finite for each $(i,j,k) \in \mU$. Also, let $\alpha, \beta, \gamma$ be formal indeterminates.  Define $f_{i,j,k} = x_{i,j,k}$ for $(i,j,k) \in \mI$ and $$f_{i,j,k} = \frac{\alpha f_{i-1,j,k}f_{i,j-1,k-1} + \beta f_{i,j-1,k}f_{i-1,j,k-1} + \gamma f_{i,j,k-1}f_{i-1,j-1,k}}{f_{i-1,j-1,k-1}}$$ for $(i,j,k) \in \mU$.  Then each $f_{i,j,k}$ is a Laurent polynomial in the $x_{i,j,k}$ with coefficients in $\Z[\alpha,\beta,\gamma]$. \end{theorem}
\pf By our usual translation and intersection argument, we may replace $\mL$ by $\mL \cap C(0,0,0)$ and define $f_{i,j,k}$ only for $(i,j,k) \in C(0,0,0)$, and it suffices to show that $f_{0,0,0}$ is a Laurent polynomial in the $x_{i,j,k}$ with coefficients in $\Z[\alpha,\beta,\gamma]$.  Working over $\Q(\sqrt{\alpha},\sqrt{\beta},\sqrt{\gamma})$, let $a_{j,k} = \sqrt{\beta\gamma/\alpha}, b_{i,k} = \sqrt{\gamma\alpha/\beta}, c_{i,j} = \sqrt{\alpha\beta/\gamma}$ for all $i, j, k$; the recurrence (\ref{recur}) then assumes the desired form.  By Theorem \ref{mainthm}, we know that $$f_{0,0,0} = \sum_G \sqrt{\beta\gamma/\alpha}^{n_a(G)} \sqrt{\gamma\alpha/\beta}^{n_b(G)} \sqrt{\alpha\beta/\gamma}^{n_c(G)} \cdot X(G)$$ (where $X(G)$ is some Laurent monomial in the variables $x_{i,j,k}$, with coefficient $1$) $$= \sum_G \alpha^\frac{n_b(G)+n_c(G)-n_a(G)}{2} \beta^\frac{n_c(G)+n_a(G)-n_b(G)}{2} \gamma^\frac{n_a(G)+n_b(G)-n_c(G)}{2} \cdot X(G).$$  The result now follows from Theorem \ref{triineq}. \eop
\par
\begin{figure}[htb] \begin{center} \includegraphics[width=3.5in]{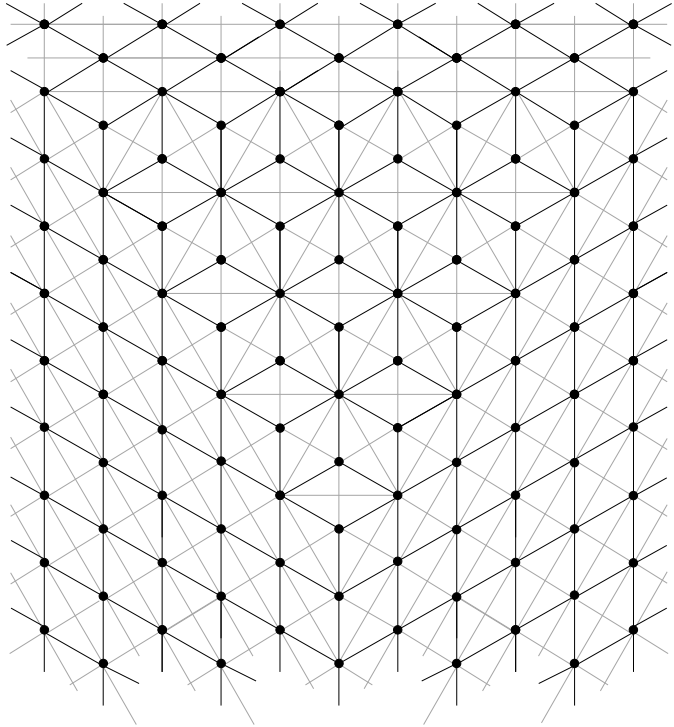} \caption{$\mG$ for Propp's standard initial conditions} \label{stdmg} \end{center} \end{figure}

 It is interesting to consider the cube recurrence in some special cases.  The question that Propp originally asked in \cite{Proppfaces} was, given the initial conditions $f_{i,j,k} = x_{i,j,k}\ (-1 \leq i+j+k \leq 1)$ and the recurrence without edge variables $$f_{i,j,k} = \frac{f_{i-1,j,k}f_{i,j-1,k-1} + f_{i,j-1,k}f_{i-1,j,k-1} + f_{i,j,k-1}f_{i-1,j-1,k}}{f_{i-1,j-1,k-1}},\quad (i+j+k > 1)$$ how to interpret combinatorially the terms of the Laurent polynomial $f_{i,j,k}$ for $i+j+k = n > 1$.  Propp observed, by setting every $x_{i',j',k'} = 1$ and using an easy induction, that there are $3^{\lfloor n^2/4\rfloor}$ such terms; the question was then what relevant objects have the property that there are $3^{\lfloor n^2/4 \rfloor}$ objects of order $n$.  We now have the means to answer this question.  By translation and $C(0,0,0)$-intersection, the problem is equivalent to describing $f_{0,0,0}$ given the initial conditions obtained from $\mL = \{(i,j,k) \in C(0,0,0)\ |\ i+j+k \leq 1-n\}$.  We get $\mI = \{(i,j,k) \in C(0,0,0)\ |\ -1-n \leq i+j+k \leq 1-n, \hbox{ or } i+j+k < -1-n \hbox{ and } \max\{i,j,k\} = 0\}$.  Thus, the terms of the polynomial correspond to groves on these initial conditions, which we call {\em standard groves of order $n$}.  As an example, a typical graph $\mG$ in which the groves embed (here for $n = 5$) is shown in Figure \ref{stdmg}.
The grove shown in Figure \ref{stdgrov} is one of the $3^{\lfloor 5^2/4\rfloor} = 729$ standard groves of order $5$.
\par
\begin{figure}[hbt] \begin{center} \includegraphics[width=6in]{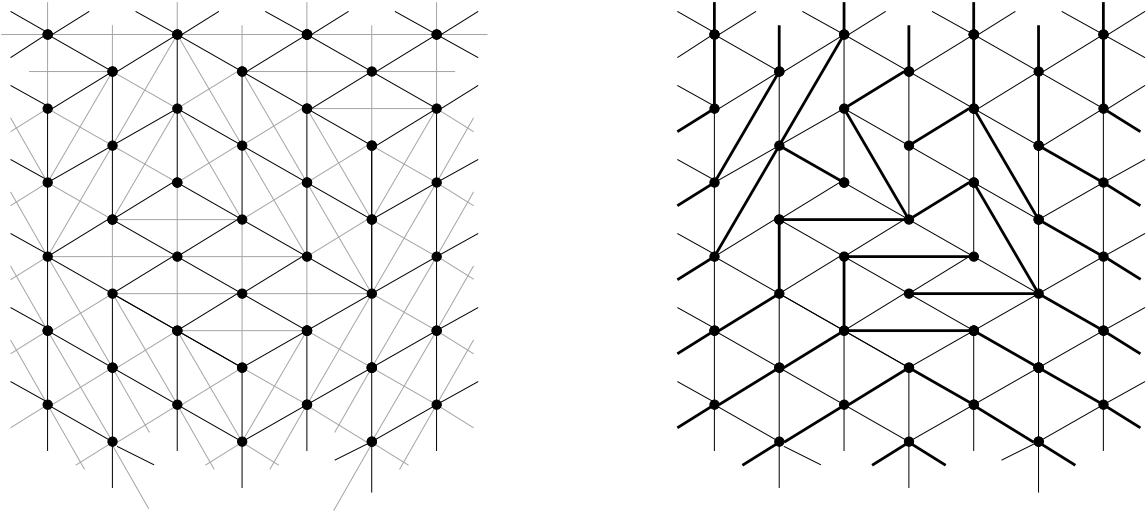} \caption{Michael Kleber's initial conditions and a sample grove} \label{klebboth} \end{center} \end{figure}

 Michael Kleber (cited in \cite{FomZel}) proposed using the initial conditions $f_{i,j,k} = x_{i,j,k}$ for $\min\{i,j,k\} = 0$, and computing $f_{i,j,k}$ (again without edge variables) for $i, j, k > 0$.  Under our usual series of manipulations, this is equivalent to finding $f_{0,0,0}$ with initial conditions induced by $$\mL = C(0,0,0) - \{(i',j',k')\ |\ i' > -i, j' > -j, k' > -k\}.$$  A typical graph $\mG$ on these initial conditions is shown in Figure \ref{klebboth} (here for $(i,j,k) = (2,3,2)$), together with a grove on these initials.
\par Another interesting special case arises in connection with the Gale-Robinson Theorem, conjectured in \cite{Gale} and proved by Fomin and Zelevinsky in \cite{FomZel} in the following form.  (The reduction of the Gale-Robinson recurrence to the cube recurrence, also used in \cite{FomZel}, was first suggested by Propp.)
\begin{theorem}[Gale-Robinson] Let $p, q, r$ be positive integers, let $n = p+q+r$, and let $\alpha, \beta, \gamma$ be formal indeterminates.  Suppose that the sequence $y_0, y_1, y_2, \ldots$ satisfies $$y_{l+n} = \frac{\alpha y_{l+p}y_{l+n-p} + \beta y_{l+q}y_{l+n-q} + \gamma y_{l+r}y_{l+n-r}}{y_l}$$ for $l \geq 0$.  Then each $y_l$ is a Laurent polynomial in the initial terms $y_0, \ldots, y_{n-1}$ with coefficients in $\Z[\alpha,\beta,\gamma]$. \end{theorem}
Indeed, for fixed $l$, we can set $\mL = \{(i,j,k) \in \Z^3\ |\ pi+qj+rk < n-l\}$, yielding $\mI = \{(i,j,k) \in \Z^3\ |\ -l \leq pi+qj+rk < n-l\}$.  After setting $x_{i,j,k} = y_{pi+qj+rk+l}\ ((i,j,k) \in \mI)$ and applying the form of the recurrence stated in Theorem \ref{fzcube}, we then obtain $f_{i,j,k} = y_{pi+qj+rk+l}$ for $(i,j,k) \in \mU$ by induction.  The desired Laurentness of $y_l = f_{0,0,0}$ is then a direct consequence of Theorem \ref{fzcube}.  Moreover, if $y_l$ is viewed as a Laurent polynomial in $y_0, \ldots, y_{n-1}, \alpha, \beta, \gamma$, then all the coefficients are nonnegative; this is immediate from Theorem \ref{mainthm} (since $y_l$ is a sum of monomials, each with coefficient $1$) and again addresses a conjecture of Fomin and Zelevinsky in \cite{FomZel}.
\par
\begin{figure}[htb] \begin{center} \includegraphics[width=3.5in]{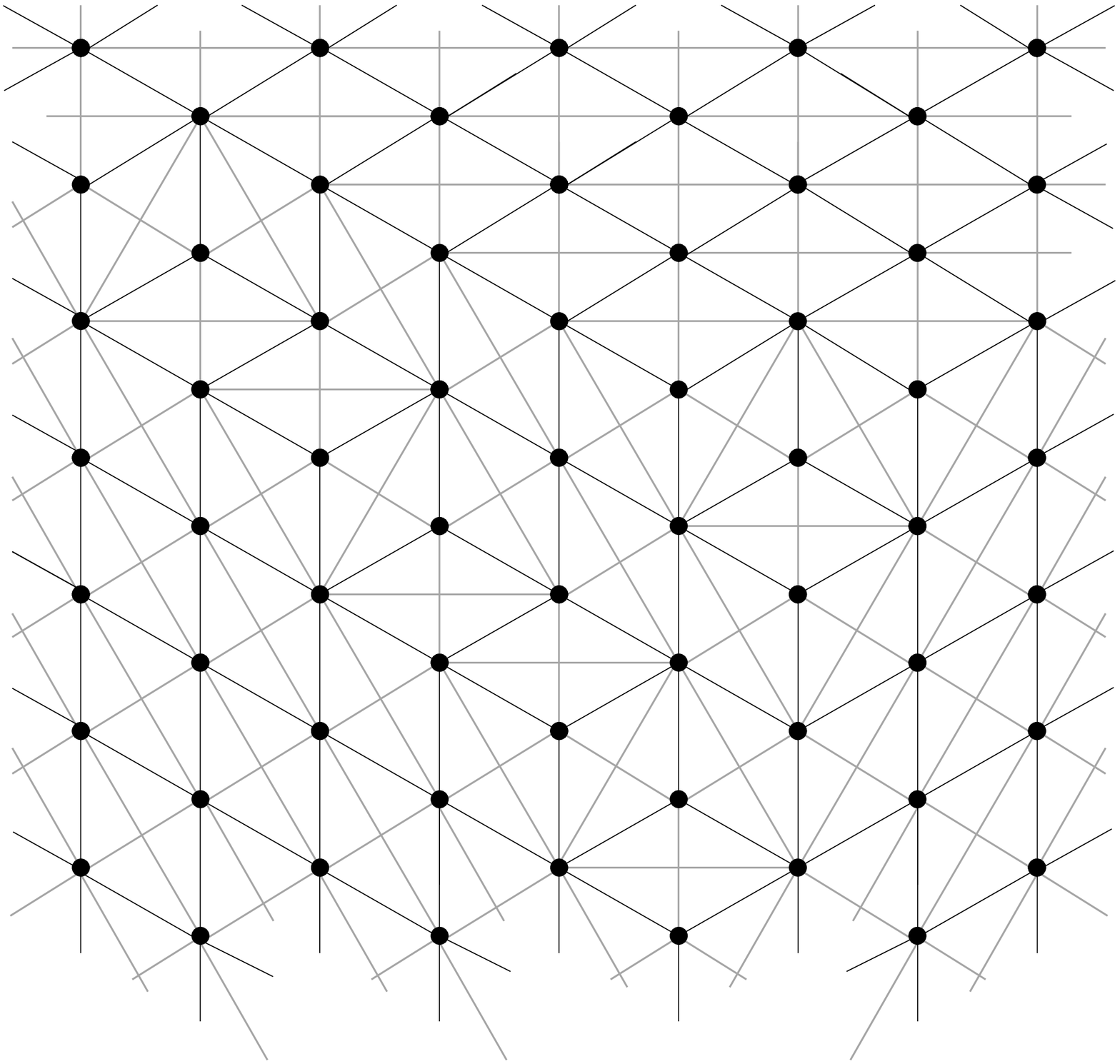} \caption{$\mG$ for a Gale-Robinson term} \label{grmg} \end{center} \end{figure}

 If we apply our usual $C(0,0,0)$-intersection, we obtain a combinatorial interpretation for the terms of any Gale-Robinson sequence in terms of groves.
In particular, the Somos-$6$ and Somos-$7$ sequences (see \cite{Gale}) can be interpreted as counting groves, resulting in a new proof that the terms of these sequences are all integers.  Figure \ref{grmg} shows the graph $\mG$ for the term $y_{11}$ of the Somos-$7$ sequence ($p = 4, q = 1, r = 2$).
\par One more specialization that merits investigation is obtained by taking arbitrary initial conditions and setting $a_{j,k} = t, b_{i,k} = t, c_{i,j} = 1/t$, as in the proof of Theorem \ref{triineq}.  Then the cube recurrence takes the form $$f_{i,j,k} = \frac{f_{i-1,j,k}f_{i,j-1,k-1} + f_{i,j-1,k}f_{i-1,j,k-1} + t^2f_{i,j,k-1}f_{i-1,j-1,k}}{f_{i-1,j-1,k-1}}.$$  In each monomial $m(G)$ coming from any polynomial $f_{i,j,k}$, the total exponent of $t$ resulting from this substitution is $n_a + n_b - n_c$, where $n_a, n_b, n_c$ are defined in terms of $G$ as specified in Theorem $\ref{triineq}$.  Since this quantity is always nonnegative, we conclude that $t$ only appears to nonnegative powers in any $f_{i,j,k}$.  Therefore, we can legitimately substitute $t = 0$, and the resulting Laurent polynomials (in the $x_{i,j,k}$) are given by the recurrence \begin{equation} \label{shiftocta} f_{i,j,k} = \frac{f_{i-1,j,k}f_{i,j-1,k-1} + f_{i,j-1,k}f_{i-1,j,k-1}}{f_{i-1,j-1,k-1}}. \end{equation}  The terms in $f_{0,0,0}$ then correspond to precisely those groves in which the total number of long edges of the forms $e_a(i,j,k)$ and $e_b(i,j,k)$ equals the number of long edges of the form $e_c(i,j,k)$.  On the other hand, by setting $g_{x,y,z} = f_{(x+y-z)/2,(x-y+z)/2,(x-y-z)/2}$ for $(x,y,z) \in \Z^3, x + y + z \equiv 0 \mod{2}$, we obtain from (\ref{shiftocta}) the recurrence $$g_{x,y,z} = \frac{g_{x-1,y-1,z}g_{x-1,y+1,z} + g_{x-1,y,z-1}g_{x-1,y,z+1}}{g_{x-2,y,z}}.$$  This is the octahedron recurrence; as shown by Speyer in \cite{Speyer}, the resulting Laurent polynomials can be interpreted as enumerating the perfect matchings of certain planar bipartite graphs (determined by $\mI$), including, for suitable initial conditions, the Aztec diamond graphs and the pine-cone graphs of \cite{BMW}.  Consequently, we have a bijection between the set of matchings of a graph (determined by $\mI$) and a particular subset of the groves on the initial conditions $\mI$.  The correspondence is not yet fully understood, but we hope to discuss the consequences in a forthcoming paper.

\section{The main proof}
\par We now turn to the proof of Theorem \ref{mainthm} and, with it, that of Lemma \ref{rhombcount}.   We maintain the assumptions $\mL \subseteq C(0,0,0)$ and $|\mU \cap C(0,0,0)| < \infty$.  Our basic technique will be induction on $|\mU \cap C(0,0,0)|$: we hold $N$ fixed and observe what happens as the initial conditions vary.  We begin with the observation that any set of initial conditions has a ``local minimum.''
\begin{lemma} \label{minimum} Suppose that $(0,0,0) \notin \mI$.  Then there exist $i, j, k \leq 0$ such that $(i-1,j,k)$, $(i,j-1,k), (i,j,k-1), (i,j-1,k-1), (i-1,j,k-1), (i-1,j-1,k), (i-1,j-1,k-1) \in \mI$ (and so $(i,j,k) \in \mU$). \end{lemma}
\pf By finiteness, we can choose $(i,j,k) \in \mU \cap C(0,0,0)$ with $i+j+k$ minimal; we claim that these values of $i, j, k$ suffice.  For example, minimality ensures that $(i-1,j,k) \in \mL$, while $(i,j+1,k+1) \notin \mL$ (otherwise we would have $(i,j,k) \in \mL$).  Therefore, $(i-1,j,k) \in \mI$.  By similar reasoning, all of the other specified points lie in $\mI$. \eop
\par Our next task will be to prove Lemma \ref{rhombcount}.  We first recall the statement.
\setcounter{slemma}{1}
\addtocounter{slemma}{-1}
\begin{slemma} Let $N$ be a cutoff for $\mI$, and let $\mJ = \{(i,j,k) \in \mI\ |\ i + j + k \geq -N-2\}$.  Then the set $\mJ$ consists of $3{N+3\choose 2} + 1$ points and contains exactly $3{N+2\choose 2}$ rhombi. \end{slemma}
\pf Fix $N$; we use induction on $|\mU \cap C(0,0,0)|$.  If this cardinality is $0$, then $\mL = C(0,0,0)$, $\mI = \{(i,j,k) \in \Z^3\ |\ \max\{i,j,k\} = 0\}$, and $\mJ = \{(i,j,k)\ |\ \max\{i,j,k\} = 0, i+j+k \geq -N-2\}$.  Each point of the form $(0,j,k), j+k \geq -N$ gives rise to a rhombus $r_a(0,j,k)$, and there are ${N+2\choose 2}$ such points.  Similarly, there are ${N+2\choose 2}$ rhombi of each of the other two types, for a total of $3{N+2\choose 2}$.  Also, it is straightforward to verify that there are $3{N+3\choose 2} + 1$ points in $\mJ$.
\par
\begin{figure}[htb] \begin{center} \includegraphics[width=4in]{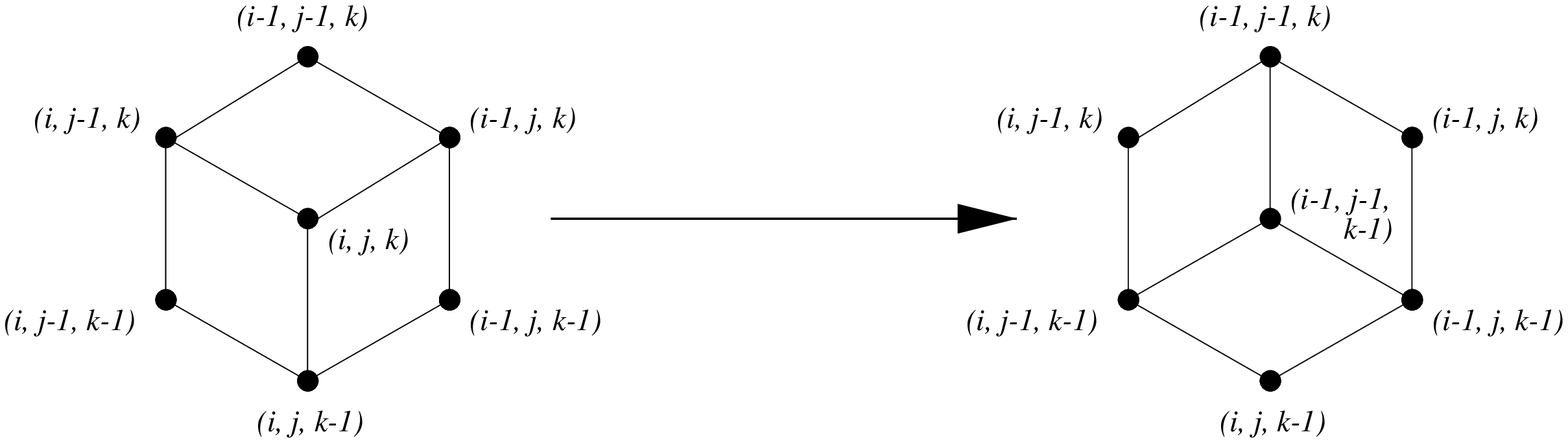} \caption{From $\mI'$ to $\mI$} \label{cubepop} \end{center} \end{figure}

 Now suppose $\mU \cap C(0,0,0) \neq \emptyset$, so that $(0,0,0) \notin \mI$.  Choose $(i,j,k)$ as given by Lemma \ref{minimum}.  Let $\mL' = \mL \cup \{(i,j,k)\}$; the lemma implies that $\mL'$ still meets the requirements we have imposed on $\mL$.  Then define $\mU', \mI', \mJ'$ by analogy with $\mU, \mI, \mJ$, so that $\mU' = \mU - \{(i,j,k)\}, \mI' = (\mI \cup \{(i,j,k)\}) - \{(i-1,j-1,k-1)\}, \mJ' = (\mJ \cup \{(i,j,k)\}) - \{(i-1,j-1,k-1)\}$.  We have $|\mU' \cap C(0,0,0)| = |\mU \cap C(0,0,0)| - 1$, and $N$ is still a cutoff for $\mI'$, so the induction hypothesis tells us that $\mJ'$ consists of $3{N+3\choose 2}+1$ points and contains $3{N+2\choose 2}$ rhombi.  However, $\mI$ is obtained from $\mI'$ (and $\mJ$ is obtained from $\mJ'$) by replacing $(i,j,k)$ with $(i-1,j-1,k-1)$.  This replacement preserves the number of points; it also preserves the number of rhombi, since it destroys the three rhombi $r_a(i,j,k), r_b(i,j,k), r_c(i,j,k)$ and creates the three rhombi $r_a(i-1,j,k), r_b(i,j-1,k), r_c(i,j,k-1)$, and one can check that $(i,j,k)$ and $(i-1,j-1,k-1)$ cannot belong to any other rhombi.  (Figure \ref{cubepop} shows the effect of this replacement in plane projection.)
\par So $\mJ$ has the same numbers of points and of rhombi as $\mJ'$, completing the induction. \eop
\par This same process also can be used to show that every set of initial conditions really does correspond to a rhombus tiling of the plane.
\par
\begin{figure}[hbt] \begin{center} \includegraphics[width=3.5in]{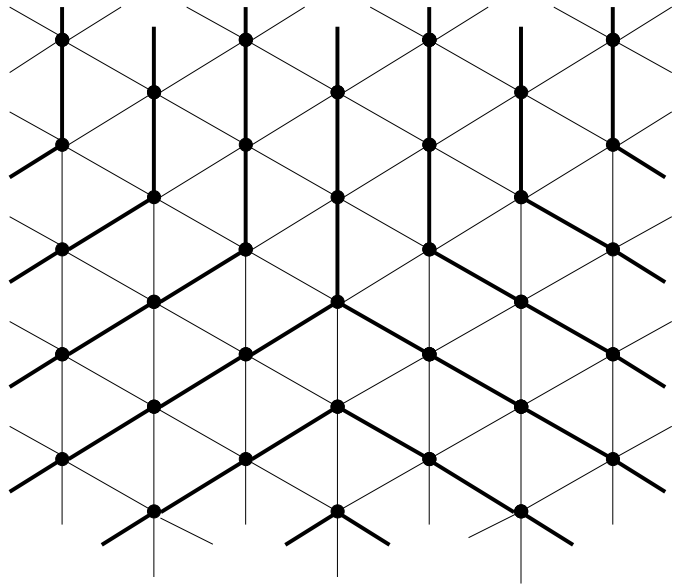} \caption{The unique grove for the base case} \label{trivgrove} \end{center} \end{figure}

 The proof of Theorem \ref{mainthm} is an application of the same technique as that of Lemma \ref{rhombcount}.  The concept is analogous to that of the ``urban renewal'' proof in \cite{Speyer}: we show that varying the initial conditions in a controlled manner is tantamount to successively applying variable substitutions in $f_{0,0,0}$, and we interpret these substitutions combinatorially.
\pfof{Theorem \ref{mainthm}} We again fix $N$ and induct on $|\mU \cap C(0,0,0)|$.  If the intersection is empty, then $\mI = \{(i,j,k)\ |\ \max\{i,j,k\} = 0\}$, and $\mG$ is as shown in Figure \ref{trivmg}.  It is straightforward to check that the graph on $\mI$ consisting of all short edges (shown in Figure \ref{trivgrove}) is a grove within radius $N$, and the corresponding monomial is $x_{0,0,0} = f_{0,0,0}$.  We claim that there are no other such groves; the base case of the induction will then follow.
\par To see this, notice that the rhombi in $\mI$ are precisely those of the forms $r_a(0,j,k)$, $r_b(i,0,k)$, and $r_c(i,j,0)\ (i, j, k \leq 0)$.  Suppose there exists a grove in which some long edge of the form $e_a(0,j,k)$ appears, and consider such an edge with $j+k$ minimal (only finitely many long edges occur).  Then we know that $j + k \geq -N$, and each of the rhombi $r_a(0,j-l-1,k-l), r_a(0,j-l,k-l-1)\ (l \geq 0)$ contributes its short edge.  In particular, we have the infinite chain of edges $$\cdots \leftrightarrow (0,j-3,k-2) \leftrightarrow (0,j-2,k-1) \leftrightarrow (0,j-1,k) \qquad$$ $$\quad \leftrightarrow (0,j,k-1) \leftrightarrow (0,j-1,k-2) \leftrightarrow (0,j-2,k-3) \leftrightarrow \cdots$$
\par However, the connectivity condition tells us, for some $l$, that $(0,j-l-1,k-l)$ should not be connected to $(0,j-l,k-l-1)$.  Thus, we have a contradiction.  So no grove (of any radius) can use any long edge of the form $e_a(0,j,k)$, and similarly the other two long edge types cannot be used.  Hence, the only possible grove is that consisting exclusively of short edges.  The base case follows.
\par Now suppose $\mU \cap C(0,0,0)$ is not empty.  Let $(i,j,k)$ be as given by Lemma \ref{minimum}, and define $\mL', \mU', \mI'$ as in the proof of Lemma \ref{rhombcount}; we again observe that $N$ remains a cutoff for $\mI'$.  Let $f'_{i',j',k'}$, for $(i',j',k') \in (\mI' \cup \mU') \cap C(0,0,0)$, be the Laurent polynomials (in $\{a_{j',k'}, b_{i',k'}, c_{i',j'}, x_{i',j',k'}\ |\ (i',j',k') \in \mI'\}$) generated by the cube recurrence from the initial conditions $\mI'$.  Because $|\mU' \cap C(0,0,0)| = |\mU \cap C(0,0,0)| - 1$, we know by induction that $f'_{0,0,0} = \sum m(G')$, where the sum is over all $\mI'$-groves $G'$ within radius $N$.  Let $g = \sum m(G)$, where the sum is over all $\mI$-groves $G$ within radius $N$.  We wish to show that $f_{0,0,0} = g$.  On the other hand, we know that $f_{0,0,0}$ is obtained from $f'_{0,0,0}$ by the variable substitution \begin{equation} \label{subst} x_{i,j,k} \leftarrow \frac{b_{i,k}c_{i,j}x_{i-1,j,k}x_{i,j-1,k-1} + c_{i,j}a_{j,k}x_{i,j-1,k}x_{i-1,j,k-1} + a_{j,k}b_{i,k}x_{i,j,k-1}x_{i-1,j-1,k}}{x_{i-1,j-1,k-1}}, \end{equation} since the right-hand side equals $f_{i,j,k}$, the left-hand side equals $f'_{i,j,k}$, and all other initial conditions are the same for $f$ as for $f'$ (and the recurrences are identical).
\par
\begin{figure}[hptb] \begin{center} \includegraphics[width=3.1in]{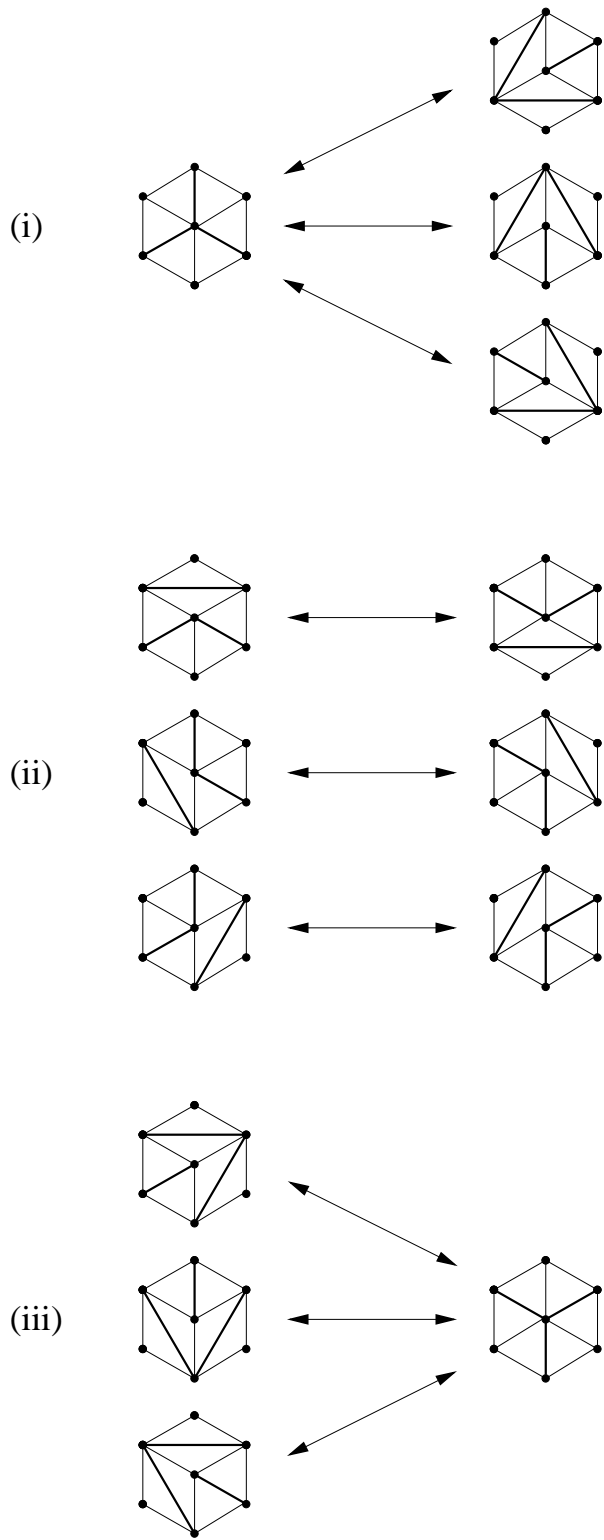} \caption{The correspondence between $\mI'$-groves (left) and $\mI$-groves (right)} \label{corresp} \end{center} \end{figure}

 We construct a correspondence between $\mI'$-groves and $\mI$-groves (within radius $N$); this correspondence is sometimes one-to-one, sometimes one-to-three, and sometimes three-to-one.  It is defined as follows: given an $\mI'$-grove, we consider the edges used in the three rhombi containing $(i,j,k)$ and replace them by any of the corresponding sets of edges in the three rhombi containing $(i-1,j-1,k-1)$, as shown in Figure \ref{corresp}; the rest of the grove is left intact.  The reverse operation --- turning $\mI$-groves into $\mI'$-groves --- is defined similarly.  In view of the connectivity constraint, we see that Figure \ref{corresp} accounts for all possible groves on either set of initial conditions (the only other possible configurations of edges would leave the vertex $(i-1,j-1,k-1)$ or $(i,j,k)$ isolated).
It is evident that our operation preserves the compactness condition, once we observe that $i + j + k > -N$.  Moreover, since the correspondence preserves all connectivity relations among vertices other than $(i,j,k)$ and $(i-1,j-1,k-1)$, it quickly follows that the connectivity condition is preserved as well, and the other conditions are trivial.  Thus, our operation does take $\mI'$-groves to $\mI$-groves, and vice versa.
\par Now let us express this correspondence algebraically.  Consider any $\mI'$-grove $G'$.  Since $(i,j,k)$ only belongs to three rhombi there, it has degree $3, 2$, or $1$ (never $0$).  If it has degree $3$, then $x_{i,j,k}$ has exponent $1$ in $m(G')$.  There are three corresponding $\mI$-groves (case (i) in Figure \ref{corresp}), which we call $G_1, G_2, G_3$, in the order in which they appear in the figure.  In $G_1$, vertex $(i-1,j-1,k-1)$ has degree $1$; vertices $(i-1,j,k)$ and $(i,j-1,k-1)$ each have degree $1$ greater than in $G'$; and new long edges $e_b(i,j-1,k)$ and $e_c(i,j,k-1)$ are used.  All other vertices and edges are the same in $G_1$ as in $G'$.  Thus, $$m(G_1)/m(G') = b_{i,k}c_{i,j}x_{i-1,j,k}x_{i,j-1,k-1}/x_{i-1,j-1,k-1}x_{i,j,k}.$$  Performing similar analyses for $G_2$ and $G_3$ (and using the fact that $x_{i,j,k}$ appears in $m(G')$ with exponent $1$), we see that $m(G_1) + m(G_2) + m(G_3)$ is obtained from $m(G')$ by the substitution (\ref{subst}).
\par Now suppose $(i,j,k)$ has degree $2$ in $G'$; thus, $x_{i,j,k}$ does not occur in $m(G')$.  From cases (ii) of Figure \ref{corresp}, we see that there is one corresponding grove $G$, and it is easy to check that $m(G) = m(G')$ --- no vertices change degrees (and $(i,j,k)$, $(i-1,j-1,k-1)$ both have degree $2$), and the one long edge that disappears is replaced by another edge represented by the same variable.  Since $x_{i,j,k}$ does not occur in $m(G')$, we may say that $m(G)$ is obtained from $m(G')$ by the substitution (\ref{subst}).
\par Finally, if $(i,j,k)$ has degree $1$, then $G'$ belongs to a triple $\{G_1, G_2, G_3\}$ of $\mI'$-groves, all of which correspond to the same $\mI$-grove $G$, as shown in case (iii) of Figure \ref{corresp}.  (We know that $G_1, G_2, G_3$ really are all $\mI$-groves because they are obtainable via our correspondence from $G$, which in turn is obtainable from our original grove $G'$.)  An analysis similar to that used in case (i) shows that $$m(G_1)/m(G) = b_{i,k}c_{i,j}x_{i-1,j,k}x_{i,j-1,k-1}/x_{i-1,j-1,k-1}x_{i,j,k}.$$  Similar computations with $G_2$ and $G_3$ give that $m(G_1) + m(G_2) + m(G_3)$ equals $$\frac{m(G)}{x_{i,j,k}} \cdot \frac{b_{i,k}c_{i,j}x_{i-1,j,k}x_{i,j-1,k-1} + c_{i,j}a_{j,k}x_{i,j-1,k}x_{i-1,j,k-1} + a_{j,k}b_{i,k}x_{i,j,k-1}x_{i-1,j-1,k}}{x_{i-1,j-1,k-1}},$$ and, using the fact that $x_{i,j,k}$ has exponent $-1$ in each of $m(G_1), m(G_2), m(G_3)$, we conclude that $m(G)$ is obtained from $m(G_1) + m(G_2) + m(G_3)$ by the same substitution (\ref{subst}).
\par Thus, summing over all $\mI'$-groves and all $\mI$-groves within radius $N$, we see that $g$ is obtained from $f'_{0,0,0}$ by applying (\ref{subst}).  Since this same substitution produces $f_{0,0,0}$ from $f'_{0,0,0}$, we have $f_{0,0,0} = g$, and the induction is complete. \eop
\par It is worth noting that a speedy combinatorial proof of Theorem \ref{acyclic} can be obtained by the same inductive process we used to prove Theorem \ref{mainthm} --- acyclicity is preserved when an $\mI'$-grove is replaced by a corresponding $\mI$-grove.  (The same is true of Theorem \ref{coeffone}.)  We provided the more geometric proof earlier in part because it elucidates better how acyclicity follows directly from the other properties of a grove.

\section{Further questions}
\par We now have one proof of Theorem \ref{mainthm}.  However, in \cite{Speyer}, Speyer gives two proofs of the analogous statement for the octahedron recurrence (in terms of perfect matchings of graphs).  One proof examines a Laurent polynomial with fixed indices under varying initial conditions, as we have done; the second proof uses fixed initial conditions and evaluates successive polynomials, making use the method of graphical condensation introduced by Eric Kuo (e.g. in \cite{Kuo}).  Briefly stated, the method operates as follows:  Suppose that we wish to determine the number of matchings of some graph $G$.  We choose some induced subgraph $H$, and we also find several pairs of subgraphs of $G$ such that, in each pair, the intersection of the two subgraphs is $H$.  When a perfect matching of $G$ and a perfect matching of $H$ are superimposed, the result is a multiset of edges of $G$.  Under suitable conditions, this multiset can then be decomposed again as a union of matchings of one of our pairs of subgraphs, with possibly some fixed set of additional edges associated with the pair (independent of the particular matchings).  If we can show that, for any such multiset, the number of decompositions into a matching of $G$ and a matching of $H$ equals the number of decompositions into matchings of one of the other pairs of subgraphs, then we obtain an expression for the number of matchings of $G$ in terms of the numbers of matchings of the subgraphs.  By associating an appropriate Laurent monomial to each matching, we can obtain a similar statement for matching-counting polynomials.
\par
\begin{figure}[hpt] \begin{center} \includegraphics[width=6in]{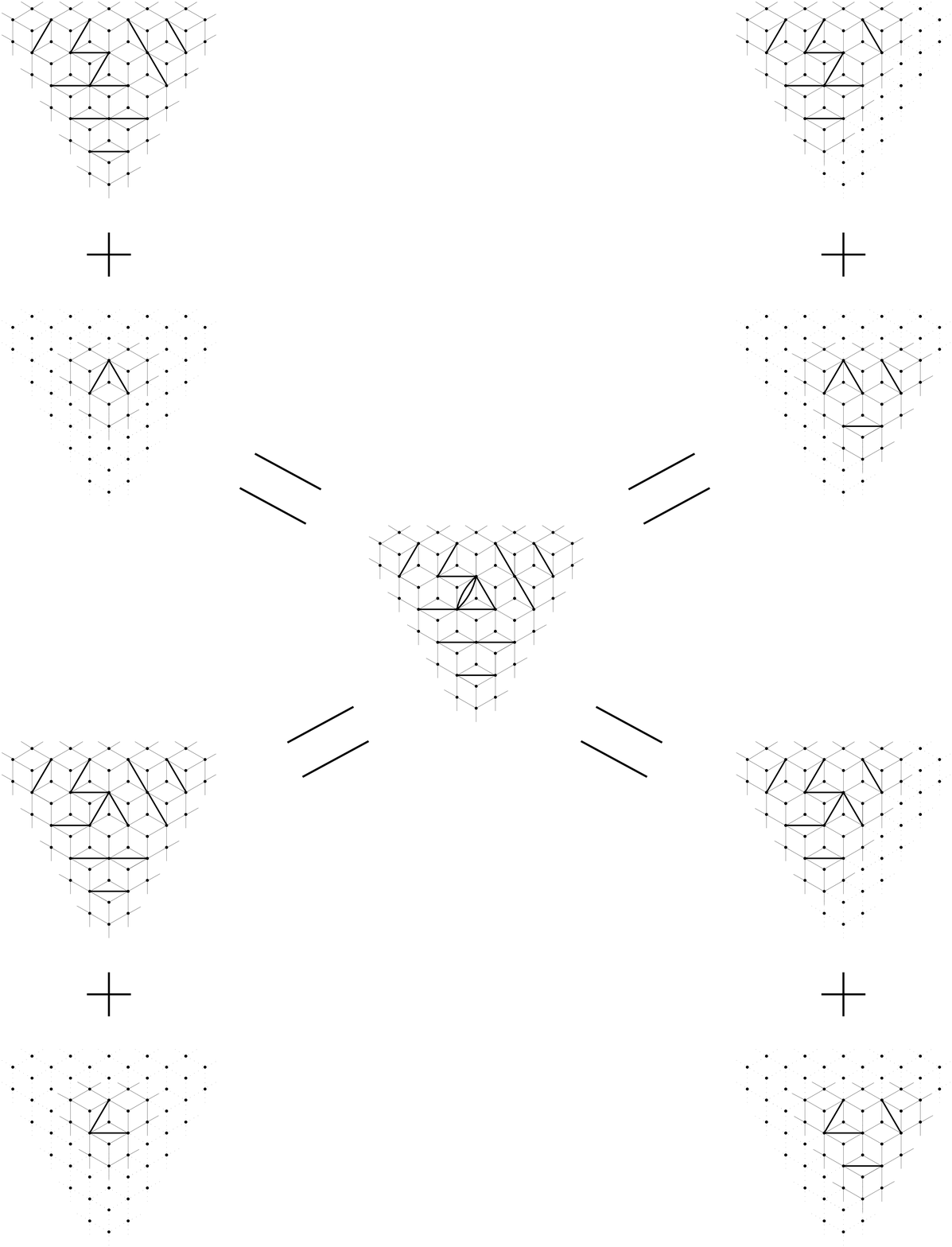} \caption{Decompositions of a doubled grove}  \label{kuogrov} \end{center} \end{figure}

 One would expect an analogous statement to hold true for groves: loosely speaking, by combining the edges of a grove $G$ with those of another grove on a different set of initial conditions (allowing for some translation), we obtain a multiset of edges; this multiset should roughly decompose again as the union of the edge sets of two groves on some other sets of initial conditions, plus some extra edges.  Indeed, because {\em every} $f_{i,j,k}$ counts groves on the initial conditions induced by a translation of $\mL \cap C(i,j,k)$, the cube recurrence (\ref{recur}) implies that the grove-condensation statement must be true.  An example is shown in Figure \ref{kuogrov}, here for the initial conditions given by $\mL = \{(i,j,k) \in C(0,0,0)\ |\ i+j+k \leq -4\}$.  For visual simplicity, the short edges are omitted; the figure shows a multiset of long edges, together with two decompositions into groves on one pair of sets of initial conditions (left) and two decompositions into groves on another pair (right).  Notice that the decompositions on the right side of the figure lack the upper-right edge and the bottom edge; these are the fixed extra edges for this pair of initials.
\par However, we have not found a proof of Theorem \ref{mainthm} along these lines.  The difficulty seems to lie in characterizing the objects obtained by combining two groves on different initial conditions.  A union of matchings of a graph and a subgraph is a relatively nice object --- it is a multiset of edges such that certain vertices each belong to exactly two edges and all other vertices belong to exactly one; consequently, it can be characterized as a vertex-disjoint union of cycles, doubled edges, paths, and isolated edges.  In contrast, no simple description of the objects obtained by combining two groves seems readily available, nor is it apparent what the procedure should be for decomposing these superimpositions back into groves.  It would be interesting to have an understanding of these gadgets, so that Theorem \ref{mainthm} could be proven by condensation.
\par It appears, however, that the machinery of Kuo condensation for groves would have to be substantially different from that known for matchings of graphs.  In the usual form of Kuo condensation, each multiset of edges can be decomposed into matchings of only one pair of subgraphs (aside from the pair $\{G,H\}$ itself).  We originally conjectured, by analogy, that each superimposition of two groves can only be decomposed back onto one of the three relevant other pairs of initial conditions.  Preliminary investigations have found this conjecture to be false.  This makes the prospect of Kuo condensation for groves all the more interesting, but deeper explorations of the topic are beyond the scope of the present paper.
\par Another open area of research concerning groves relates to the connection between the cube and octahedron recurrences discussed above.  Does the correspondence between perfect matchings and groves tell us anything about matchings --- or about groves?  In particular, the known Kuo condensation algorithm for matchings may shed light on the analogue for groves.
\par
\begin{figure}[bthp] \begin{center} \includegraphics[width=6in]{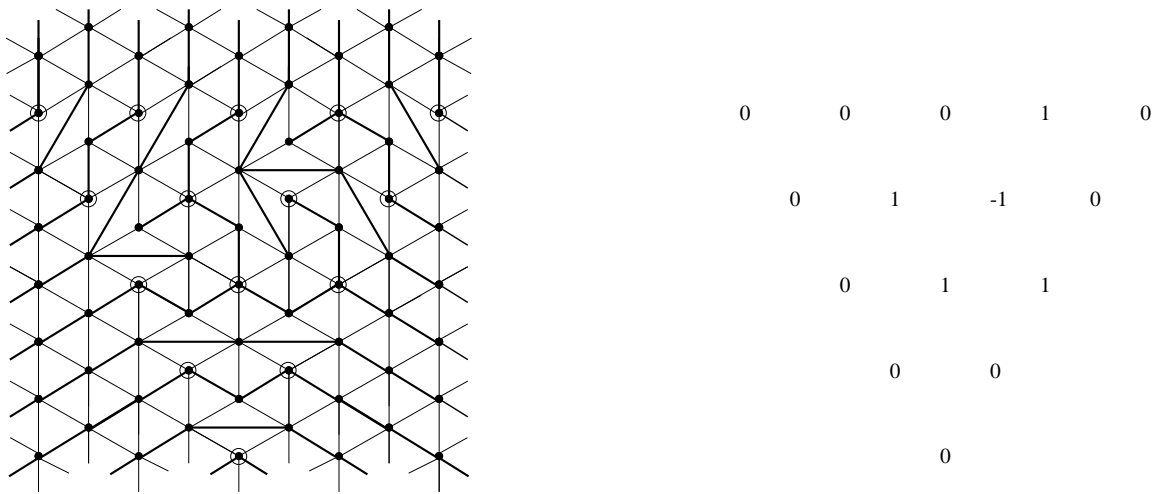} \caption{A standard grove of order $4$ and the corresponding alternating-sign triangle} \label{astri} \end{center} \end{figure}

 One more interesting question relates to the form of the cube recurrence originally proposed by Propp.  The polynomials generated by the octahedron recurrence, with standard initial conditions $\{(i,j,k) \in \Z^3\ |\ i+j+k \equiv 0 \mod{2}; k \in \{-1,0\}\}$, can be seen not only as enumerating matchings of graphs but also as enumerating compatible pairs of alternating-sign matrices: in each monomial, the exponents with which the variables $x_{i,j,k}$ appear correspond to the entries of the matrices.  (See \cite{RobRum} for details.)  In particular, if we set $x_{i,j,k} = 1$ whenever $k = -1$, then the terms of the resulting polynomial (in the variables $x_{i,j,0}$) correspond precisely to single alternating-sign matrices.  As Propp observed in \cite{Proppfaces}, something analogous seems to be taking place with the cube recurrence for initial conditions $\mI = \{(i,j,k) \in \Z^3\ |\ -1 \leq i+j+k \leq 1\}$.  Specifically, fix $x_{i,j,k} = 1$ for $i+j+k = -1, 0$, and also set each edge variable equal to $1$.  Then any $f_{i,j,k}$ is a Laurent polynomial in $\{x_{i',j',k'}\ |\ i'+j'+k' = 1; i' \leq i, j' \leq j, k' \leq k\}$.  For any monomial, the exponents of these $x_{i',j',k'}$, which must be $-1, 0$, or $1$ (this is easy to check, since, after translation and intersection, each such point only belongs to three rhombi), form a sort of ``alternating-sign triangle.''  An example is shown in Figure \ref{astri}: at left is the grove for a particular monomial; at right is the alternating-sign triangle for the same monomial, formed by the exponents of the $x_{i',j',k'}$ corresponding to the circled vertices.
The question therefore arises as to whether there exists some nice description for the set of alternating-sign triangles of a given size, analogous to that for alternating-sign matrices.  The prospect of developing a theory of alternating-sign triangles is a tantalizing one.
\par The analogy with matchings of graphs thus provides a wellspring of inspiration for further questions about groves.  We hope that the results herein established will constitute the foundation for fruitful future study.

\section{Acknowledgments}
\par We would like to thank the other members of the REACH (Research Experiences in Algebraic Combinatorics at Harvard) group for assistance in the early stages of investigating the cube recurrence and for comments on drafts of this paper.  We would also like to thank the National Security Agency and the National Science Foundation for their financial support of our research, as well as the mathematics departments of Harvard University and the University of Wisconsin-Madison for administrative assistance.  Most of all, we thank Jim Propp, the director of REACH, for introducing us to the cube recurrence and for guiding and directing our work on the subject.

\end{document}